\documentclass[11pt]{article}

\usepackage{amsfonts,amssymb,amsmath}
\usepackage[english]{babel}
 \def\botcaption#1#2{\medskip\centerline{{\scshape #1.}\kern8pt
 {\rm #2}}\bigskip}

 \makeatletter
 \@addtoreset{equation}{subsection}
 \makeatother

 \makeatletter
 \@addtoreset{enunciato}{section}
 \makeatother

 \newcounter{enunciato}[subsection]

 \newtheorem{ittheorem}{Theorem}
 \newtheorem{itlemma}{Lemma}
 \newtheorem{itproposition}{Proposition}
 \newtheorem{itdefinition}{Definition}
 \newtheorem{itremark}{Remark}
 \newtheorem{itclaim}{Claim}

 \newenvironment{theorem}{\addtocounter{enunciato}{1}
 \begin{ittheorem}}{\end{ittheorem}}

 \newenvironment{lemma}{\addtocounter{enunciato}{1}
 \begin{itlemma}}{\end{itlemma}}

 \newenvironment{proposition}{\addtocounter{enunciato}{1}
 \begin{itproposition}}{\end{itproposition}}

 \newenvironment{definition}{\addtocounter{enunciato}{1}
 \begin{itdefinition}}{\end{itdefinition}}

 \newenvironment{remark}{\addtocounter{enunciato}{1}
 \begin{itremark}}{\end{itremark}}

 \newenvironment{claim}{\addtocounter{enunciato}{1}
 \begin{itclaim}}{\end{itclaim}}

 \newenvironment{proof}{\noindent {\bf Proof.\,}
 }{\hspace*{\fill}$\square$\medskip}

 \newcommand{\be}[1]{\begin{equation}\label{#1}}
 \newcommand{\ee}{\end{equation}}

 \newcommand{\bl}[1]{\begin{lemma}\label{#1}}
 \newcommand{\el}{\end{lemma}}

 \newcommand{\br}[1]{\begin{remark}\label{#1}}
 \newcommand{\er}{\end{remark}}

 \newcommand{\bt}[1]{\begin{theorem}\label{#1}}
 \newcommand{\et}{\end{theorem}}

 \newcommand{\bd}[1]{\begin{definition}\label{#1}}
 \newcommand{\ed}{\end{definition}}

 \newcommand{\bcl}[1]{\begin{claim}\label{#1}}
 \newcommand{\ecl}{\end{claim}}

 \newcommand{\bp}[1]{\begin{proposition}\label{#1}}
 \newcommand{\ep}{\end{proposition}}

 \newcommand{\bc}[1]{\begin{corollary}\label{#1}}
 \newcommand{\ec}{\end{corollary}}

 \newcommand{\bpr}{\begin{proof}}
 \newcommand{\epr}{\end{proof}}

 \newcommand{\bi}{\begin{itemize}}
 \newcommand{\ei}{\end{itemize}}

 \newcommand{\ben}{\begin{enumerate}}
 \newcommand{\een}{\end{enumerate}}

 \def\botcaption#1#2{\medskip\centerline{{\scshape #1.}\kern8pt
 {\rm #2}}\bigskip}

 \def \ba {\begin{array}}
 \def \ea {\end{array}}

 \def \Z {{\mathbb Z}}
 \def \R {{\mathbb R}}
 
 \def \N {{\mathbb N}}
 
 \def \da {\downarrow}
 \def \ua {\uparrow}

 \def \cW {{\mathcal W}}
 \def \cD {{\mathcal D}}
 \def \cL {{\mathcal L}}
 \def \cI {{\mathcal I}}
 \def \cR {{\mathcal R}}
 \def \cA {{\mathcal A}}
 \def \cE {{\mathcal E}}
 \def \cP {{\mathcal P}}

 \def \DOM {{\hbox{\footnotesize\rm DOM}}}
 \def \CONE {{\hbox{\footnotesize\rm CONE}}}

\begin{document}

\title{Localization transition for a copolymer in an emulsion}

\author{F.\ den Hollander
\footnote{Mathematical Institute, Leiden University, P.O.\ Box 9512, 
2300 RA Leiden, The Netherlands}\,
\footnote{EURANDOM, P.O.\ Box 513, 5600 MB Eindhoven, The Netherlands}\\
S.G.\ Whittington
\footnote{Department of Chemistry, University of Toronto, Toronto, 
Canada M5S 3H6}
}

\maketitle

\begin{center}
{\sl Dedicated to Ya.G.\ Sinai on the occasion of his 70th birthday}
\end{center}

\begin{abstract}
In this paper we study a two-dimensional directed self-avoiding walk model of a 
random copolymer in a random emulsion. The polymer is a random concatenation of 
monomers of two types, $A$ and $B$, each occurring with density $\frac{1}{2}$. 
The emulsion is a random mixture of liquids of two types, $A$ and $B$, organised
in large square blocks occurring with density $p$ and $1-p$, respectively, where
$p \in (0,1)$. The polymer in the emulsion has an energy that is minus $\alpha$ 
times the number of $AA$-matches minus $\beta$ times the number of $BB$-matches, 
where $\alpha,\beta\in\R$ are interaction parameters. Symmetry considerations 
show that without loss of generality we may restrict to the cone $\{(\alpha,\beta)
\in\R^2\colon\,\alpha\geq |\beta|\}$.

We derive a variational expression for the quenched free energy per monomer in 
the limit as the length $n$ of the polymer tends to infinity and the blocks in 
the emulsion have size $L_n$ such that $L_n \to \infty$ and $L_n/n \to 0$. To 
make the model mathematically tractable, we assume that the polymer can only 
enter and exit a pair of neighbouring blocks at diagonally opposite corners. 
Although this is an unphysical restriction, it turns out that the model exhibits 
rich and physically relevant behaviour. 

Let $p_c \approx 0.64$ be the critical probability for directed bond percolation 
on the square lattice. We show that for $p \geq p_c$ the free energy has a 
phase transition along \emph{one} curve in the cone, which turns out to be 
\emph{independent} of $p$. At this curve, there is a transition from a phase 
where the polymer is fully $A$-delocalized (i.e., it spends almost all of its 
time deep inside the $A$-blocks) to a phase where the polymer is partially 
$AB$-localized (i.e., it spends a positive fraction of its time near those 
interfaces where it diagonally crosses the $A$-block rather than the $B$-block). 
We show that for $p<p_c$ the free energy has a phase transition along \emph{two} 
curves in the cone, both of which turn out to \emph{depend} on $p$. At the 
\emph{first} curve there is a transition from a phase where the polymer is 
fully $A,B$-delocalized (i.e., it spends almost all of its time deep inside 
the $A$-blocks and the $B$-blocks) to a partially $BA$-localized phase, while 
at the \emph{second} curve there is a transition from a partially $BA$-localized 
phase to a phase where both partial $BA$-localization and partial $AB$-localization 
occur simultaneously. 

We derive a number of qualitative properties of the critical curves. The 
supercritical curve is non-decreasing and concave with a finite horizontal 
asymptote. Remarkably, the first subcritical curve does not share these
properties and does not converge to the supercritical curve as $p \ua p_c$.
Rather, the second subcritical curve converges to the supercritical curve
as $p \da 0$.

\vskip 0.5truecm
\noindent
\emph{AMS} 2000 \emph{subject classifications.} 60F10, 60K37, 82B27.
\end{abstract}


\section{Introduction and main results}
\label{S1}

\subsection{Background}
\label{S1.1}

(Linear) copolymers are polymer chains consisting of two or more types of monomer. 
Random copolymers are copolymers where the order of the monomers along the 
polymer chain is determined by a random process. In any particular chain, 
the sequence of monomers once determined is fixed, so a random copolymer is an 
example of a \emph{quenched random system}. In this paper we will be concerned 
with copolymers consisting of two types of monomer, labelled $A$ and $B$. We 
write $\omega_i\in\{A,B\}$ to denote the type of the $i$-th monomer, and 
$\omega = \{\omega_1,\omega_2,\cdots\}$ to denote the full order along the chain, 
which is truncated at $\omega_n$ when the polymer has length $n$.  We will only 
consider the case where the random variables $\omega_i$ are independent and 
identically distributed (i.i.d.), assuming the values $A$ and $B$ with probability
$\frac12$ each. In principle, the properties of the polymer depend on $\omega$, 
and we write $P^\omega_n$ for the value of a property $P$ when the polymer has 
length $n$ and order $\omega$. If $\lim_{n\to\infty} P^\omega_n$ exists $\omega$-a.s.\ 
and is non-random, then we say that the property $P$ is \emph{self-averaging}.

Several different physical situations are of interest. For instance, if the
monomer-monomer interactions differ for pairs $AA$, $BB$ and $AB$, then we
may investigate the effect of the randomness on the \emph{collapse transition},
where the polymer collapses from a random coil to a ball-like structure as 
the temperature decreases or the solvent quality varies. Alternatively, if the 
two types of monomer interact differently with an impenetrable surface, then 
we may investigate the effect of the randomness on the \emph{adsorption transition},
where the polymer adsorbs onto the surface as the temperature decreases or the 
surface quality varies. There are interesting questions about how the location 
of the collapse transition or the adsorption transition, and the values of 
associated critical exponents, depend on the parameters controlling the randomness. 
Many of these questions remain unresolved. For background and references, the reader 
is referred to Orlandini \emph{et al} \cite{OTW99}, \cite{OTW00}, Janse van Rensburg 
\emph{et al} \cite{JvROTW01}, Brazhnyi and Stepanow \cite{BS02}, Whittington
\cite{W02}, and Soteros and Whittington \cite{SW04}.

The problem that we will consider here is the \emph{localization transition} of 
a random copolymer near an interface. Suppose that we have two immiscible liquids 
and that it is energetically favourable for monomers of one type to be in one 
liquid and for monomers of the other type to be in the other liquid. At high 
temperatures the polymer will delocalize into one of the liquids in order to 
maximise its entropy, while at low temperatures energetic effects will dominate 
and the polymer will localize close to the interface between the two liquids in 
order to be able to place more than half of its monomers in their preferred liquid. 
In the limit as $n\to\infty$, we may expect a phase transition. A typical example 
here would be an oil-water interface and a copolymer with hydrophobic and hydrophilic 
monomers.

Given such a physical situation, the polymer can be modelled in a variety
of ways, e.g.\ as a random walk or as a self-avoiding walk, either directed
or undirected. Such examples have been investigated for the situation where 
\emph{the interface is flat and infinite}. In addition, there is some flexibility 
in the details of the Hamiltonian that is chosen to model the interactions.
 
A simple model with a single interface was proposed and analysed by Garel 
\emph{et al} \cite{GHLO89}, with a Hamiltonian that depends on temperature 
and interaction bias. A first mathematical treatment of this model was given 
by Sinai \cite{S93} and by Grosberg \emph{et al} \cite{GIN94}, in the absence 
of interaction bias, for a directed random walk version of the model. For this 
version, Bolthausen and den Hollander \cite{BdH97} proved that the quenched free 
energy is non-analytic along a critical curve in the plane of inverse temperature 
vs.\ interaction bias, and derived several qualitative properties of this curve,
among which upper and lower bounds. Albeverio and Zhou \cite{AZ96}, for the 
unbiased case, and Biskup and den Hollander \cite{BdH99}, for the biased case, 
extended this work by deriving path properties of the model, in particular, 
ergodic limits along the interface and exponential tightness perpendicular to 
the interface in the localized phase, as well as zero limiting frequency of hits 
of the interface in the interior of the delocalized phase. The latter result 
was recently strengthened by Giacomin and Toninelli \cite{GT05}, who showed 
that in the interior of the delocalized phase the number of times the path 
intersects the interface grows at most logarithmically with its length. The 
conjecture is that the number of intersections is actually bounded. In Giacomin 
and Toninelli \cite{GT05b} it was proved that the free energy is infinitely 
differentiable inside the localized phase. Thus, there is no phase transition of 
finite order anywhere off the critical curve. Morover, in Giacomin and Toninelli 
\cite{GT05a} it was proved that the free energy is twice differentiable across 
the critical curve, i.e., the phase transition is at least of second order.

Maritan \emph{et al} \cite{MRT99} considered both random walk and self-avoiding 
walk models and derived rigorous bounds on the free energy, under an assumption 
on the asymptotics of a certain class of self-avoiding walks. Martin \emph{et al} 
\cite{MCW00} proved the existence of a localization transition for a self-avoiding 
walk model and obtained qualitative results about the shape of the phase transition 
curve. These results were extended and improved by Madras and Whittington \cite{MW03}, 
who also gave a rigorous version of the result of Maritan \emph{et al} \cite{MRT99} 
for the self-avoiding walk model. Orlandini \emph{et al} \cite{ORW02} derived 
rigorous bounds on the critical curve for the directed random walk model, 
while Causo and Whittington \cite{CW03} and James, Soteros and Whittington 
\cite{JSW03} obtained sharp numerical estimates for the self-avoiding walk model.

An interesting recent development concerns the slope of the critical curve 
in the limit of small inverse temperature and interaction bias in the directed 
random walk version of the model. In Bolthausen and den Hollander \cite{BdH97} 
it was proved that this slope exists, is strictly positive and is at most $1$, the 
latter being a corollary of an upper bound on the full critical curve. Garel 
\emph{et al} \cite{GHLO89} had earlier hinted at the possibility that this slope be 
1, a viewpoint that was taken up by Trovato and Maritan \cite{TM99}. However, 
Stepanow \emph{et al} \cite{SSE98} conjectured the slope to be $\frac23$, based 
on replica symmetry arguments. Monthus \cite{M00}, using a general renormalization 
scheme, conjectured a simple explicit formula for the full critical curve, 
which indeed has slope $\frac23$ in the limit of small inverse temperature and 
interaction bias. Based on this work, Bodineau and Giacomin \cite{BG04pr} proved 
that this formula is a lower bound for the critical curve, so that we now know 
that the slope is at least $\frac23$. Numerical work by Garel and Monthus 
\cite{GM05} and Caravenna, Giacomin and Gubinelli \cite{CGG05} indicates that 
the upper and lower bounds on the critical curve are not sharp, nor are the 
bounds $1$ and $\frac23$ for the slope. So far all attempts to improve these 
bounds have failed. The slope seems to be close to $0.82$.

The reason why the above issue is of interest is that, while the full shape of 
the critical curve is model-dependent, the slope in the limit of small 
inverse temperature and interaction bias is believed to be insensitive to the 
details of the model.

The goal of the present paper is to study a model where \emph{the interface has a 
more complex geometry}. A first attempt in this direction was made by den Hollander
and W\"uthrich \cite{dHW04}, where \emph{an infinite array of parallel flat infinite
interfaces} was considered and the average hopping time between interfaces was 
computed for a directed random walk model. In the present paper we investigate 
the situation in which the lattice is divided into large blocks, and each block is 
independently labelled $A$ or $B$ with probability $p$ and $1-p$, respectively,
i.e., \emph{the interface has a percolation type structure}. This is a primitive 
model of an emulsion (e.g.\ oil dispersed as droplets in water as the dispersing 
medium). As before, the copolymer consists of a random concatenation of monomers 
of type $A$ and $B$. It is energetically favourable for monomers of type $A$ to 
be in the $A$-blocks and for monomers of type $B$ to be in the $B$-blocks of the
emulsion. \emph{Under the restriction that the polymer can only enter and exit a
pair of neighbouring blocks at diagonally opposite corners}, we show that there 
is a phase transition between a phase where the polymer is \emph{fully delocalized} 
away from the interfaces between the two types of blocks and a phase where the 
polymer is \emph{partially localized} near the interfaces. It turns out that the 
critical curve \emph{does not} depend on $p$ in the supercritical percolation 
regime, but \emph{does} depend on $p$ in the subcritical percolation regime. In 
the latter regime, a \emph{second} critical curve appears separating \emph{two} 
partially localized phases. 

Our paper is organised as follows. In the rest of Section \ref{S1} we define 
the model, formulate our main theorems, discuss these theorems, and formulate
some open problems. Section \ref{S2} contains some preparatory results about 
path entropies and free energies per pair of neighbouring blocks. In Sections 
\ref{S3} and \ref{S4} we provide the proofs of the main theorems, focussing
on the free energy, respectively, the critical curves.

\subsection{The model}
\label{S1.2}

Each positive integer is randomly labelled $A$ or $B$, with probability $\frac{1}{2}$ 
each, independently for different integers. The resulting labelling is denoted by
\be{bondlabel}
\omega = \{\omega_i \colon\, i \in \N\} \in \{A,B\}^\N.
\ee
Fix $p \in (0,1)$ and $L_n \in \N$. Partition $\R^2$ into square blocks of size $L_n$:
\be{blocks}
\R^2 = \bigcup_{x \in \Z^2} \Lambda_{L_n}(x), \qquad
\Lambda_{L_n}(x) = xL_n + (0,L_n]^2.
\ee
(Note that the blocks contain their north and east side but not their south and west 
side.) Each block is randomly labelled $A$ or $B$, with probability $p$, respectively, 
$1-p$, independently for different blocks. The resulting labelling is denoted by
\be{blocklabel}
\Omega = \{\Omega(x) \colon\, x \in \Z^2\} \in \{A,B\}^{\Z^2}.
\ee

Consider the set of $n$-step \emph{directed self-avoiding paths} starting at the 
origin and being allowed to move \emph{upwards, downwards and to the right}. Let 
$\cW_{n,L_n}$ be the subset of those paths that enter blocks at a corner, exit 
blocks at one of the two corners \emph{diagonally opposite} the one where it entered, 
and in between \emph{stay confined} to the two blocks that are seen when entering. 
In other words, after the path reaches a site $xL_n$, it must make a step to the 
right, it must subsequently stay confined to the pair of blocks labelled $x$ and 
$x+(0,-1)$, and it must exit this pair of blocks either at site $xL_n+(L_n,L_n)$ 
or at site $xL_n+(L_n,-L_n)$ (see Figure 1). This restriction is put in to make 
the model mathematically tractable.


\vskip 1.2truecm

\setlength{\unitlength}{0.25cm}

\begin{picture}(10,10)(-18,-5)

  {\thicklines
   \qbezier(0,6)(3,6)(6,6)
   \qbezier(6,0)(6,3)(6,6)
   \qbezier(0,0)(3,0)(6,0)
   \qbezier(6,-6)(6,-3)(6,0)
  }
  \qbezier[40](0,0)(0,3)(0,6)
  \qbezier[40](0,0)(0,-3)(0,-6)
  \qbezier[40](0,-6)(3,-6)(6,-6)
  \put(6,6){\circle*{0.5}}
  \put(6,-6){\circle*{0.5}}
  \put(0,0){\circle*{0.5}}
  \put(-3.3,0){$xL_n$}
  \put(4,7){$xL_n+(L_n,L_n)$}
  \put(4,-7.8){$xL_n+(L_n,-L_n)$}  

  \put(-15,-10){\small 
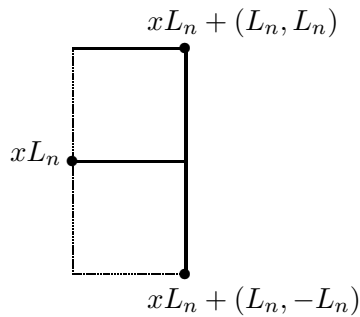
               Fig.\ 1. Two neighbouring blocks. The dots are the sites of
               entrance and exit.
               \normalsize}
  \put(-15,-11.5){\small
               The drawn lines are part of the blocks, the dashed lines are not.
               \normalsize}
\end{picture}

\vskip 2truecm


Given $\omega,\Omega$ and $n$, with each path $\pi \in \cW_{n,L_n}$ we associate 
an energy given by the Hamiltonian
\be{Hamiltonian}
H_{n,L_n}^{\omega,\Omega}(\pi) 
= - \sum_{i=1}^n \Big(\alpha 1\{\omega_i=\Omega^{L_n}_{\pi_i}=A\}
+ \beta 1\{\omega_i=\Omega^{L_n}_{\pi_i}=B\}\Big),
\ee 
where $\pi_i$ denotes the $i$-th step of the path and $\Omega^{L_n}_{\pi_i}$ 
denotes the label of the block that step $\pi_i$ lies in. What this Hamiltonian 
does is count the number of $AA$-matches and $BB$-matches and assign them energy 
$-\alpha$ and $-\beta$, respectively, where $\alpha,\beta\in\R$. Note that 
\emph{the interaction is assigned to bonds rather than to sites}: we identify 
the monomers with the steps of the path. 

For $\alpha,\beta>0$, the above definitions are to be interpreted as follows: 
$\omega$ plays the role of the random monomer types, with $A$ denoting hydrophobic 
and $B$ denoting hydrophilic; $\Omega$ plays the role of the random emulsion, with 
$A$ denoting oil and $B$ denoting water; $n$ is the number of monomers; the 
Hamiltonian assigns negative energy to matches of affinities between polymer 
and emulsion, with $\alpha$ and $\beta$ the interaction strengths (it assigns 
zero energy to mismatches).

Given $\omega,\Omega$ and $n$, we define the \emph{quenched free energy per step} as
\be{fedef}
\begin{aligned}
f_{n,L_n}^{\omega,\Omega}
&= \frac{1}{n} \log Z_{n,L_n}^{\omega,\Omega},\\
Z_{n,L_n}^{\omega,\Omega}
&= \sum\limits_{\pi\in\cW_{n,L_n}} \exp\left[-H_{n,L_n}^{\omega,\Omega}(\pi)\right].
\end{aligned}
\ee
We are interested in the limit $n \to \infty$ subject to the restriction
\be{Ln}
L_n \to \infty \qquad \mbox{ and } \qquad L_n/n \to 0.
\ee
This is a \emph{coarse-graining} limit where the path spends a long time in each 
single block yet visits many blocks. Throughout the paper we assume that this 
restriction is in force, which is necessary to make the model mathematically 
tractable. It will turn out that the free energy does not depend on the choice of 
the sequence $(L_n)_{n\in\N}$.

\subsection{Free energy}
\label{S1.3}

Theorem \ref{feiden} below says that the quenched free energy per step is 
self-averaging and can be expressed as a variational problem involving the 
free energies of the polymer in each of the four possible pairs of adjacent 
blocks it may encounter and the frequencies at which the polymer visits each 
of these pairs of blocks on the coarse-grained block scale. To formulate this 
theorem we need some more definitions.

First, for $L\in\N$ and $a \geq 2$ (with $aL$ integer), let $\cW_{aL,L}$ 
denote the set of $aL$-step directed self-avoiding paths starting at $(0,0)$, 
ending at $(L,L)$, and in between not leaving the two adjacent blocks of size 
$L$ labelled $(0,0)$ and $(-1,0)$.


\vskip 1.2truecm

\setlength{\unitlength}{0.25cm}

\begin{picture}(10,10)(-18,-5)

  {\thicklines
   \qbezier(0,6)(3,6)(6,6)
   \qbezier(6,0)(6,3)(6,6)
   \qbezier(0,0)(3,0)(6,0)
   \qbezier(6,-6)(6,-3)(6,0)
  }
  \qbezier[40](0,0)(0,3)(0,6)
  \qbezier[40](0,0)(0,-3)(0,-6)
  \qbezier[40](0,-6)(3,-6)(6,-6)
  \put(6,6){\circle*{0.5}}
  \put(6,-6){\circle*{0.5}}
  \put(0,0){\circle*{0.5}}
  \put(-4,0){$(0,0)$}
  \put(4,7){$(L,L)$}
  \put(4,-7.8){$(L,-L)$}  
  
  \qbezier[45](0,0)(3,3)(6,6)
  {\thicklines
   \qbezier(2.5,3)(2.75,3)(3,3)
   \qbezier(3,2.5)(3,2.75)(3,3)
  }

  \put(-11,-10){\small 
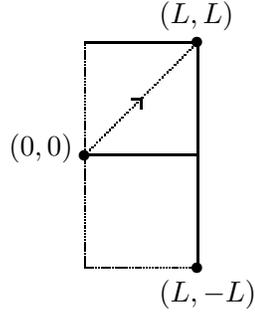
               Fig.\ 2. Two neighbouring blocks. The dashed line with arrow
               \normalsize}
  \put(-11,-11.5){\small
               indicates that the coarse-grained path makes a step diagonally 
               upwards.
               \normalsize}
\end{picture}

\vskip 2truecm


\noindent
For $k,l\in\{A,B\}$, let
\be{feklaL}
\begin{aligned}
\psi^\omega_{kl}(aL,L) &= \frac{1}{aL} \log Z^{\omega}_{aL,L},\\
Z^{\omega}_{aL,L} &= \sum_{\pi\in\cW_{aL,L}}
\exp\big[-H^{\omega,\Omega}_{aL,L}(\pi)\big]
\hbox{ when } \Omega(0,0)=k \hbox{ and } \Omega(0,-1)=l,\\
\end{aligned} 
\ee
denote the free energy per step in a $kl$-block when the number of steps inside
the block is $a$ times the size of the block. Let
\be{fekl}
\lim_{L\to\infty} \psi^\omega_{kl}(aL,L) = \psi_{kl}(a) 
= \psi_{kl}(\alpha,\beta;a).
\ee
Note here that $k$ labels the type of the block that is diagonally crossed, while
$l$ labels the type of the block that appears as its neighbour at the starting
corner (see Fig.\ 2). We will prove in Section \ref{S2.2} that the limit exists 
$\omega$-a.s.\ and is non-random. It will turn out that $\psi_{AA}$ and $\psi_{BB}$ 
take on a simple form, whereas $\psi_{AB}$ and $\psi_{BA}$ do not.

Second, let $\cW$ denote the class of all \emph{coarse-grained paths} $\Pi=
\{\Pi_j\colon\,j\in\N\}$ that step diagonally from corner to corner (see Fig.\ 3, 
where each dashed line with arrow denotes a single step of $\Pi$). For $n\in\N$, 
$\Pi\in\cW$ and $k,l\in\{A,B\}$, let 
\be{fracdef}
\begin{array}{ll}
\rho^\Omega_{kl}(\Pi,n) = \frac{1}{n}
\sum_{j=1}^n 1\,\{&\hbox{$\Pi_j$ diagonally crosses a $k$-block 
in $\Omega$ that has an $l$-block}\\
&\hbox{in $\Omega$ appearing as its neighbour at the 
starting corner}\,\,\,\}.
\end{array}
\ee
Abbreviate
\be{rhoklp}
\rho^\Omega(\Pi,n) = \left(\rho^\Omega_{kl}(\Pi,n)\right)_{k,l\in\{A,B\}}, 
\ee
which is a $2 \times 2$ matrix with nonnegative elements that sum up to 1.
Let $\cR^\Omega(\Pi)$ denote the set of all limits points of the sequence
$\{\rho^\Omega(\Pi,n)\colon\,n\in\N\}$, and put
\be{Ldef}
\cR^\Omega = \mbox{the closure of the set } \bigcup_{\Pi\in\cW} \cR^\Omega(\Pi).
\ee
Clearly, $\cR^\Omega$ exists for all $\Omega$. Moreover, since $\Omega$ has 
a trivial sigma-field at infinity (i.e., all events not depending on finitely
many coordinates of $\Omega$ have probability 0 or 1) and $\cR^\Omega$ is measurable 
with respect to this sigma-field, we have
\be{Rdef}
\cR^\Omega = \cR(p) \qquad \Omega-a.s.
\ee
for some \emph{non-random closed} set $\cR(p)$. This set, which depends on 
the parameter $p$ controlling $\Omega$, will be analysed in Section \ref{S3.2}. 
It is the set of all possible limit points of the frequencies at which the 
four pairs of adjacent blocks can be seen along an infinite coarse-grained path.


\vskip 2.7truecm

\setlength{\unitlength}{0.3cm}

\begin{picture}(10,10)(-10,0)

  {\thicklines
   \qbezier(0,16)(8,16)(16,16)
   \qbezier(16,0)(16,8)(16,16)
   \qbezier(0,0)(8,0)(16,0)
   \qbezier(0,0)(0,8)(0,16)
   \qbezier(4,0)(4,8)(4,16)
   \qbezier(8,0)(8,8)(8,16)
   \qbezier(12,0)(12,8)(12,16)
   \qbezier(0,4)(8,4)(16,4)
   \qbezier(0,8)(8,8)(16,8)
   \qbezier(0,12)(8,12)(16,12)
  }
  \qbezier[30](0,4)(2,6)(4,8)
  \qbezier[30](4,8)(6,10)(8,12)
  \qbezier[30](8,12)(10,10)(12,8)
  \qbezier[30](12,8)(14,6)(16,4)
  {\thicklines
   \qbezier(1.7,6)(1.85,6)(2,6)
   \qbezier(2,5.7)(2,5.85)(2,6)
   \qbezier(5.7,10)(5.85,10)(6,10)
   \qbezier(6,9.7)(6,9.85)(6,10)
   \qbezier(9.9,9.8)(10.05,9.8)(10.2,9.8)
   \qbezier(10.2,10.1)(10.2,9.95)(10.2,9.8)
   \qbezier(13.9,5.8)(14.05,5.8)(14.2,5.8)
   \qbezier(14.2,6.1)(14.2,5.95)(14.2,5.8)
  }

  \put(1,14.5){$A$}
  \put(5,14.5){$B$}
  \put(9,14.5){$A$}
  \put(13,14.5){$A$}
  \put(1,10.5){$B$}
  \put(5,10.5){$A$}
  \put(9,10.5){$B$}
  \put(13,10.5){$B$}
  
  \put(1,6.5){$B$}
  \put(5,6.5){$A$}
  \put(9,6.5){$A$}
  \put(13,6.5){$A$}
  \put(1,2.5){$B$}
  \put(5,2.5){$B$}
  \put(9,2.5){$A$}
  \put(13,2.5){$B$}
  
  \put(0.05,4){\circle*{.6}}
  \put(0.05,12){\circle*{.6}}
  \put(4.05,0){\circle*{.6}}
  \put(4.05,8){\circle*{.6}}
  \put(4.05,16){\circle*{.6}}
  \put(8.05,4){\circle*{.6}}
  \put(8.05,12){\circle*{.6}}
  \put(12.05,0){\circle*{.6}}
  \put(12.05,8){\circle*{.6}}
  \put(12.05,16){\circle*{.6}}
  \put(16.05,4){\circle*{.6}}
  \put(16.05,12){\circle*{.6}}

  \put(-9,-2.5){\small 
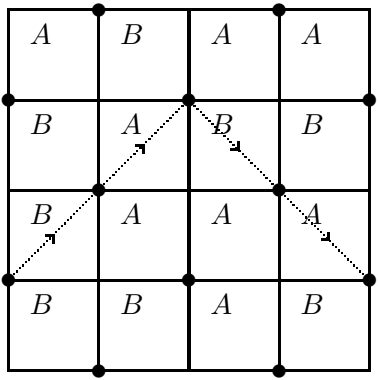
               Fig.\ 3. $\Pi$ sampling $\Omega$. 
               The dashed lines with arrows indicate the steps of $\Pi$.}
  
\end{picture}

\vskip 1.3truecm


Let $\cA$ be the set of $2 \times 2$ matrices whose elements are $\geq 2$. The
starting point of our paper is the following representation of the free energy.

\bt{feiden}
(i) For all $(\alpha,\beta)\in\R^2$ and $p \in (0,1)$,
\be{sa}
\lim_{n\to\infty} f_{n,L_n}^{\omega,\Omega} = f = f(\alpha,\beta;p)
\ee
exists $\omega,\Omega$-a.s., is finite and non-random, and is given by 
\be{fevar}
f = \sup_{(a_{kl}) \in \cA}\, \sup_{(\rho_{kl}) \in \cR(p)}
\frac{\sum_{k,l} \rho_{kl} a_{kl} \psi_{kl}(a_{kl})}
{\sum_{k,l} \rho_{kl} a_{kl}}.
\ee
(ii) The function $(\alpha,\beta)\mapsto f(\alpha,\beta;p)$ is convex on $\R^2$
for all $p \in (0,1)$.\\
(iii) The function $p\mapsto f(\alpha,\beta;p)$ is continuous on $(0,1)$ for
all $(\alpha,\beta)\in\R^2$.\\
(iv) For all $(\alpha,\beta)\in\R^2$ and $p \in (0,1)$,
\be{symms}
\begin{aligned}
f(\alpha,\beta;p) &= f(\beta,\alpha;1-p),\\
f(\alpha,\beta;p) &= \frac12(\alpha+\beta) + f(-\beta,-\alpha;p).
\end{aligned}
\ee
\et 

\noindent
Theorem \ref{feiden}(i), which will be proved in Section \ref{S3.1}, says that the 
limiting free energy per step is \emph{self-averaging in both $\omega$ and $\Omega$}, 
and equals the average of the limiting free energies per step associated with the four 
pairs of adjacent blocks, weighted and optimised according to the frequencies at which 
these four pairs are visited by the coarse-grained path and the fraction of time spent 
in each of them by the path. Assumption (\ref{Ln}) is crucial, since it provides the
separation of the path scale and the block scale, thereby separating the self-averaging
in $\omega$ and $\Omega$. Theorem \ref{feiden}(ii), which will be proved in Section 
\ref{S3.1} also, is standard. Theorem \ref{feiden}(iii) is a consequence of the
fact that $p\mapsto\cR(p)$ is continuous in the Hausdorff metric, which will be proved
in Section \ref{S3.2}. Theorem \ref{feiden}(iv) is immediate from (\ref{Hamiltonian})
upon interchanging the two monomer types and/or the two block types.

In view of Theorem \ref{feiden}(iv), we may \emph{without loss of generality restrict
to the cone}
\be{OCTdef}
\CONE = \{(\alpha,\beta)\in\R^2\colon\,\alpha \geq |\beta|\}.
\ee
The upper half of the cone is the physically most relevant part, but we will see that
also the lower half of the cone is of interest. Note that \emph{$AA$-matches are favored 
over $BB$-matches}. This will be crucial throughout the paper. 

The behaviour of $f$ as a function of $(\alpha,\beta)$ is different for $p \geq p_c$ 
and $p < p_c$, where $p_c \approx 0.64$ is the critical percolation density for directed 
bond percolation on the square lattice. The reason is that the coarse-grained paths
$\Pi$, which determine the set $\cR(p)$, sample $\Omega$ just like paths in \emph{directed 
bond percolation} on the square lattice rotated by 45 degrees sample the percolation 
configuration (see Fig.\ 3).

\subsection{Supercritical case $p \geq p_c$}
\label{S1.4}

The entropy per step of the walk in a single block, subject to (\ref{Ln}), is
\be{kappadef}
\kappa = \lim_{n\to\infty} \frac{1}{n} \log |\cW_{n,L_n}|.
\ee
In Section \ref{S2.1} we will see that $\kappa = \frac{1}{2} \log 5$. This number 
is special to our model.

Our first theorem identifies the two phases, which turn out \emph{not to depend 
on} $p$.

\bt{fesupcr}
Let $p \geq p_c$. Then $f(\alpha,\beta;p)=f(\alpha,\beta)$, and $(\alpha,\beta)\mapsto 
f(\alpha,\beta)$ is non-analytic along the curve in $\CONE$ separating the two regions
\be{DL}
\begin{aligned}
\cD &= \hbox{ delocalized regime} 
    &= \left\{(\alpha,\beta)\in\CONE \colon f(\alpha,\beta) 
     = \frac{1}{2}\alpha + \kappa\right\},\\
\cL &= \hbox{ localized regime} 
    &= \left\{(\alpha,\beta)\in\CONE \colon f(\alpha,\beta) 
     > \frac{1}{2}\alpha + \kappa\right\}.
\end{aligned}
\ee
\et

\noindent
The intuition behind Theorem \ref{fesupcr}, which will be proved in Section \ref{S4.1.1}, 
is as follows. The $A$-blocks (almost) percolate. Therefore the polymer has the option 
of moving to the (incipient) infinite cluster of $A$-blocks and staying in that infinite 
cluster forever, thus seeing only $AA$-blocks. In doing so, it loses an entropy of at 
most $o(n/L_n) = o(n)$, it gains an energy $\frac{1}{2}\alpha n+o(n)$ (because only half 
of its monomers are matched), and it gains an entropy $\kappa n+o(n)$. Alternatively, the 
path has the option of following the boundary of the infinite cluster, at least part of 
the time, during which it sees $AB$-blocks and (when $\beta\geq 0$) gains more energy 
by matching more than half of its monomers. Consequently,
\be{feineq}
f(\alpha,\beta) \geq \frac{1}{2}\alpha + \kappa.
\ee 
The boundary between the two regimes in (\ref{DL}) corresponds to the crossover where 
one option takes over from the other.

Our second theorem gives an explicit classification of the two phases in terms of
the free energies $\psi_{AA}$ and $\psi_{AB}$. 

\bt{phtriden}
Let $p \geq p_c$. Then 
\be{DLdef}
\begin{aligned}
\cD &= \{(\alpha,\beta)\in\CONE \colon\, S_{AB}=S_{AA}\},\\
\cL &= \{(\alpha,\beta)\in\CONE \colon\, S_{AB}>S_{AA}\},
\end{aligned}
\ee
where
\be{Sdef}
S_{kl}=S_{kl}(\alpha,\beta)=\sup_{a\geq 2} \psi_{kl}(\alpha,\beta;a).
\ee
\et

\noindent
We have $S_{AB} \geq S_{AA}$ for all $(\alpha,\beta)$, because in an $AB$-block the path 
may spend all of its time in the half that is $A$, in which case it is not aware of the 
presence of the half that is $B$ (see Fig.\ 4). Theorem \ref{phtriden}, which will be 
proved in Section \ref{S4.1.1} also, says that the critical curve marks those parameter 
values where $=$ changes to $>$.


\vskip 0.7truecm

\setlength{\unitlength}{0.3cm}

\begin{picture}(10,10)(-10,-4)

  {\thicklines
   \qbezier(0,6)(3,6)(6,6)
   \qbezier(6,0)(6,3)(6,6)
   \qbezier(0,0)(3,0)(6,0)
   \qbezier(6,-6)(6,-3)(6,0)
  }
  \qbezier[40](0,0)(0,3)(0,6)
  \qbezier[40](0,0)(0,-3)(0,-6)
  \qbezier[40](0,-6)(3,-6)(6,-6)
  
  {\thicklines
   \qbezier(2.7,3)(2.85,3)(3,3)
   \qbezier(3,2.7)(3,2.85)(3,3)
   \qbezier(13.15,2.9)(13.3,2.95)(13.45,3)
   \qbezier(13.45,2.7)(13.45,2.85)(13.45,3)
  }

  \put(1.5,4.5){$A$}
  \put(1.5,-1.5){$B$}
  \qbezier[60](0,0)(3,3)(6,6)
  
  {\thicklines
   \qbezier(9,6)(12,6)(15,6)
   \qbezier(15,0)(15,3)(15,6)
   \qbezier(9,0)(12,0)(15,0)
   \qbezier(15,-6)(15,-3)(15,0)
  }
  \qbezier[40](9,0)(9,3)(9,6)
  \qbezier[40](9,0)(9,-3)(9,-6)
  \qbezier[40](9,-6)(12,-6)(15,-6)
  
  \put(10.5,4.5){$A$}
  \put(10.5,-1.5){$B$}
  \qbezier[50](12,0.2)(13.5,3.1)(15,6)
  \qbezier[30](9,.2)(10.5,.2)(12,.2)

  \put(0,0){\circle*{.5}} 
  \put(6,6){\circle*{.5}}
  \put(6,-6){\circle*{.5}}
  \put(9,0){\circle*{.5}}
  \put(15,6){\circle*{.5}}
  \put(15,-6){\circle*{.5}}

  \put(-8,-8){\small 
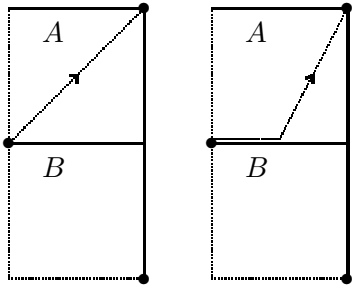
               Fig.\ 4. Two possible strategies inside an $AB$-block: The path
               can either move
               \normalsize}
  \put(-8,-9.3){\small
               straight across or move along the interface for awhile
               and then move across.
               \normalsize}
  \put(-8,-10.6){\small
                 Both strategies correspond to a coarse-grained step diagonally
                 upwards.
                 \normalsize}

\end{picture}

\vskip 2.5truecm


Our third theorem gives the qualitative properties of the critical curve
separating $\cD$ and $\cL$ (see Fig.\ 5).

\bt{phtrcurve}
Let $p \geq p_c$.\\
(i) For every $\alpha \geq 0$ there exists a $\beta_c(\alpha) \in [0,\alpha]$ 
such that the copolymer is
\be{betac}
\ba{lll}
&\mbox{delocalized} &\mbox{if }\, -\alpha \leq \beta \leq \beta_c(\alpha),\\
&\mbox{localized}   &\mbox{if }\, \beta_c(\alpha) < \beta \leq \alpha. 
\ea
\ee
(ii) The function $\alpha\mapsto\beta_c(\alpha)$ is continuous, non-decreasing and 
concave on $[0,\infty)$.\\ 
(iii) There exists an $\alpha^* \in (0,\infty)$ such that
\be{phtrlims}
\ba{lll}
&\beta_c(\alpha) = \alpha &\mbox{if }\,\alpha \leq \alpha^*,\\ 
&\beta_c(\alpha) < \alpha &\mbox{if }\,\alpha > \alpha^*.
\ea
\ee
Moreover,
\be{slopedis}
\lim_{\alpha\da\alpha^*} \frac{\alpha-\beta_c(\alpha)}{\alpha-\alpha^*} \in [0,1).
\ee
(iv) There exists a $\beta^* \in [\alpha^*,\infty)$ such that
\be{phtrlimsext}
\lim_{\alpha \to \infty} \beta_c(\alpha) = \beta^*.
\ee
\et


\setlength{\unitlength}{0.45cm}

\begin{picture}(12,12)(-6,-2)

  \put(0,0){\line(12,0){12}}
  \put(0,0){\line(0,8){8}}
  \put(0,0){\line(0,-3){3}}
  \put(0,0){\line(-4,0){4}}
  {\thicklines
   \qbezier(2,2)(4,3.2)(9,4.5)
   \thicklines
   \qbezier(0,0)(1,1)(2,2)
  }
  \qbezier[25](-.4,.5)(.3,1.2)(2,2)
  \qbezier[60](0,5.5)(5,5.5)(10,5.5)
  \qbezier[60](2,2)(4.5,4.5)(7,7)
  \qbezier[15](2,0)(2,1)(2,2)
  \qbezier[70](4,-4)(0,0)(-4,4)
  \put(-.7,-1){$0$}  
  \put(12.5,-0.2){$\alpha$}
  \put(-0.1,8.5){$\beta$}
  \put(2,2){\circle*{.3}}
  \put(1.7,-1){$\alpha^*$}
  \put(-1.4,5.4){$\beta^*$}
  \put(10.5,4.5){$\beta_c(\alpha)$}
  \put(5.6,4.3){$\cal L$}
  \put(6,2.3){$\cal D$}

  \put(-3,-5.5){\small 
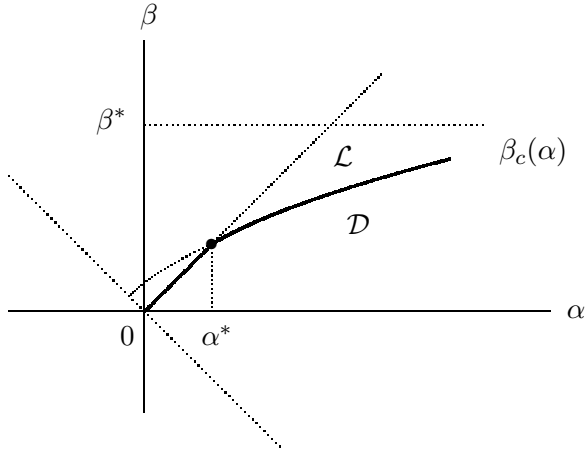
               Fig.\ 5. Qualitative picture of $\alpha\mapsto\beta_c(\alpha)$
               for $p \geq p_c$. The curved dotted line
               \normalsize}
   \put(-3,-6.4){\small
               is the analytic continuation outside $\CONE$ of the part off 
               the diagonal. 
               \normalsize}
\end{picture}

\vskip 2.5truecm


\noindent
It is clear from (\ref{DL}) that the part off the diagonal is a critical line. 
We will see in Section \ref{S4.2.3} that also the part on the diagonal is a 
critical line. Theorem \ref{phtrcurve}, which will be proved in Section 
\ref{S4.1.2}, says that the critical curve follows the diagonal for 
$\alpha\in [0,\alpha^*]$, moves off the diagonal at $\alpha=\alpha^*$ with a 
slope discontinuity, and has a finite asymptote for large $\alpha$. The concavity 
of the curve implies that it is strictly increasing as long as it is below the 
asymptote. We are not able to exclude that the curve hits the asymptote, nor 
that it follows the diagonal all the way up to the asymptote, but we expect 
this not to happen. We will see in Section \ref{S4.1.2} that the curved dotted 
line crosses the vertical axis at $(0,\alpha_0)$ with $\alpha_0 \approx 0.125$. 
We have no numerical values for $\alpha^*$ and $\beta^*$. We will show in Section 
\ref{S4.1.2} that $\beta^* \in [\log 2,8\log 3)$. Clearly, $\alpha^* \in 
[\alpha_0,\beta^*]$.

To prove Theorem \ref{phtrcurve}, we will \emph{reformulate} the criterion 
$S_{AB}>S_{AA}$ in terms of a criterion for the free energy of a model with a 
\emph{single linear interface}. This reformulation, which will be given in 
Section \ref{S2.3}, is crucial in allowing us to get a handle on the critical 
curve in Fig.\ 5.

We will see in Section \ref{S4.1.1} that $\cD$ corresponds to the situation where
the polymer is \emph{fully $A$-delocalized} (i.e., it spends almost all of its 
time away from the interface deep inside the $A$-blocks), while $\cL$ corresponds 
to the situation where the polymer is \emph{partially $AB$-delocalized} (i.e., it 
spends a positive fraction of its time near those interfaces where it diagonally 
crosses the $A$-block rather than the $B$-block).

\subsection{Subcritical case $p < p_c$}
\label{S1.5}

In the subcritical percolation regime, the analogue of the critical curve in Fig.\ 5
turns out \emph{to depend on} $p$ and to be much more difficult to characterise. 
We begin with some definitions.

Let
\be{rho*def}
\rho^*(p) = \sup_{(\rho_{kl})\in\cR(p)} [\,\rho_{AA}+\rho_{AB}\,].
\ee
This is the maximal frequency of $A$-blocks crossed by an infinite coarse-grained path
(recall (\ref{fracdef}--\ref{Rdef})). The graph of $p\mapsto\rho^*(p)$ looks like:


\vskip 1truecm

\setlength{\unitlength}{0.45cm}

\begin{picture}(8,8)(-7,-1)

  \put(0,0){\line(9,0){9}}
  \put(0,0){\line(0,7){7}}

  {\thicklines
  \qbezier(5,6)(6.25,6)(7.5,6)
  \qbezier(0,0)(1,5)(5,6)
  }

   \qbezier[40](5,0)(5,3)(5,6)
   \qbezier[40](0,6)(2.5,6)(5,6)
   \qbezier[40](7.5,0)(7.5,3)(7.5,6)
   \qbezier[60](0,0)(3.75,3)(7.5,6)
   \put(-.6,-.8){$0$}
   \put(4.8,-.8){$p_c$}
   \put(5,6){\circle*{.35}}
   \put(9.5,-.3){$p$}
   \put(-.5,7.5){$\rho^*(p)$}
   \put(-.7,5.8){$1$}
   \put(7.3,-.8){$1$}
  
  \put(-2,-2.5){\small 
               Fig.\ 6. Qualitative picture of $p\mapsto\rho^*(p)$.
               \normalsize}
  
\end{picture}

\vskip 1truecm


\noindent
Further details will be given in Section \ref{S3.2}. 

For $x,y\geq 2$, let $u(x)=u(\alpha;x)$ and $v(y)=v(\beta;y)$ be defined by
\be{fgdef}
\begin{aligned}
xu(x) &= \frac12\alpha x + \log 2 + \frac12 x\log x - \frac12(x-2)\log(x-2),\\
yv(y) &= \frac12\beta y  + \log 2 + \frac12 y\log y - \frac12(y-2)\log(y-2).
\end{aligned}
\ee
For $\rho\in(0,1)$, let
\be{Fdef}
F(\rho) =F(\alpha,\beta;\rho)
= \sup_{x,y\geq 2} \frac{\rho xu(x) + (1-\rho) yv(y)}{\rho x+ (1-\rho) y}.
\ee
This variational formula will be analysed in Section \ref{S2.5}. There we will
see that $(\alpha,\beta) \mapsto F(\alpha,\beta;\rho)$ is analytic on $\R^2$ 
for all $\rho\in (0,1)$.

The following is the analogue of Theorem \ref{fesupcr}, and will be proved in
Section \ref{S4.2.1}.

\bt{fesubcr}
Let $p<p_c$. Then $(\alpha,\beta)\mapsto f(\alpha,\beta;p)$ is non-analytic
along the curve in $\CONE$ separating the two regions
\be{DL*}
\begin{aligned}
\cD &= \hbox{ delocalized regime} 
    &= \left\{(\alpha,\beta)\in\CONE \colon f(\alpha,\beta;p) 
     = F(\alpha,\beta;\rho^*(p))\right\},\\
\cL &= \hbox{ localized regime} 
    &= \left\{(\alpha,\beta)\in\CONE \colon f(\alpha,\beta;p) 
     > F(\alpha,\beta;\rho^*(p))\right\}.
\end{aligned}
\ee
\et

\noindent
The intuition behind Theorem \ref{fesubcr} is as follows. We will see in Section
\ref{S2.2.1} that $\psi_{AA}(a) = u(a)$ and $\psi_{BB}(a) = v(a)$. In the delocalized
regime, the polymer stays away from the $AB$-interface. For the free energy this
means that no difference is felt between $\psi_{AB},\psi_{AA}$ and between $\psi_{BA},
\psi_{BB}$. Therefore in this regime the variational formula in (\ref{fevar}) 
effectively reduces to
\be{fevarred}
f = \sup_{(a_{kl})\in\cA}\,\sup_{(\rho_{kl})\in\cR(p)}
\frac{\rho_A a_{AA}\psi_{AA}(a_{AA}) + \rho_B a_{BB}\psi_{BB}(a_{BB})} 
{\rho_A a_{AA} + \rho_B a_{BB}},
\ee  
where $\rho_A=\rho_{AA}+\rho_{AB}$ and $\rho_B=\rho_{BA}+\rho_{BB}$ are the 
frequencies at which the polymer diagonally traverses $A$-blocks and $B$-blocks, 
while $a_{AA}$ and $a_{BB}$ are the respective times spent inside these blocks. 
The first supremum amounts to optimising over $a_{AA},a_{BB}\geq 2$. Since 
$AA$-matches are preferred over $BB$-matches, implying $\psi_{AA}\geq\psi_{BB}$, 
the second supremum is taken at the largest possible value of $\rho_A=1-\rho_B$ 
in $\cR(p)$, which is $\rho^*(p)$. Hence, putting $a_{AA}=x$ and $a_{BB}=y$, we 
get $f=F(\alpha,\beta;\rho^*(p))$. In the localized regime, on the other hand, 
the polymer spends part of its time near $AB$-interfaces or $BA$-interfaces, in 
which case a difference is felt between $\psi_{AB},\psi_{AA}$ and/or between 
$\psi_{BA},\psi_{BB}$, and the free energy is larger. In Section \ref{S4.2.1} 
we will make the above intuition rigorous.

Comparing the first lines of (\ref{DL}) and (\ref{DL*}), we see that the free energy 
in the supercritical delocalized regime is a function of $\alpha$ only and has a 
\emph{simple linear form}, whereas the free energy in the subcritical delocalized 
regime is a function of $\alpha,\beta,\rho^*(p)$ and has a form that is \emph{rather 
more complicated}. For $\rho=1$, (\ref{fgdef}--\ref{Fdef}) yield $F(\alpha,\beta;1) 
= \sup_{x\geq 2} u(x) = u(\frac52) = \frac12\alpha + \frac 12\log 5$. This explains 
the connection between (\ref{DL}) and (\ref{DL*}). 

\medskip
The following is the analogue of Theorem \ref{phtriden}, and will be proved in
Section \ref{S4.2.1}.

\bt{phtriden*}
Let $p<p_c$. Then
\be{DLdef*}
\begin{aligned}
\cD &= \{(\alpha,\beta)\in\CONE\colon\,\psi_{AB}(\bar x)=\psi_{AA}(\bar x)
\mbox{ and } \psi_{BA}(\bar y)=\psi_{BB}(\bar y)\},\\
\cL &= \{(\alpha,\beta)\in\CONE\colon\,\psi_{AB}(\bar x)>\psi_{AA}(\bar x)
\mbox{ or } \psi_{BA}(\bar y)>\psi_{BB}(\bar y)\},\\
\end{aligned}
\ee
where $\bar x = \bar x(\alpha,\beta;\rho^*(p))$ and $\bar y = \bar y(\alpha,\beta;\rho^*(p))$
are the unique maximisers of $F(\alpha,\beta;\rho^*(p))$, i.e., of the variational formula in
{\rm (\ref{Fdef})} for $\rho=\rho^*(p)$.
\et

\noindent
Theorem \ref{phtriden*} says that the crossover into the localized regime occurs when the 
difference between $\psi_{AB},\psi_{AA}$ or between $\psi_{BA},\psi_{BB}$ is felt at the 
minimisers of the variational formula for the delocalized regime. 

Comparing (\ref{DLdef}) and (\ref{DLdef*}), we see that the crossover into the 
supercritical localization regime occurs when the maxima of $\psi_{AB},\psi_{AA}$
separate, whereas the crossover into the subcritical localization regime occurs
when $\psi_{AB},\psi_{AA}$ or $\psi_{BA},\psi_{BB}$ separate at specific locations, 
which themselves depend on $\alpha,\beta,\rho^*(p)$.

We will see in Section \ref{S4.2.2} that $\psi_{BA}(\bar y)=\psi_{BB}(\bar y)$ implies
$\psi_{AB}(\bar x)=\psi_{AA}(\bar x)$. Hence, (\ref{DLdef*}) in fact reduces to
\be{DLdef*red}
\begin{aligned}
\cD &= \{(\alpha,\beta)\in\CONE\colon\,\psi_{BA}(\bar y)=\psi_{BB}(\bar y)\},\\
\cL &= \{(\alpha,\beta)\in\CONE\colon\,\psi_{BA}(\bar y)>\psi_{BB}(\bar y)\}.\\
\end{aligned}
\ee
This is to be interpreted as saying that, when the critical curve is crossed 
from $\cD$ to $\cL$, localization occurs in the $BA$-blocks rather than in the 
$AB$-blocks. The intuitive explanation is as follows. In the delocalized phase
the polymer spends positive fractions of its time in the $A$-blocks and in the 
$B$-blocks (the $A$-blocks do not percolate). Because $AA$-matches are preferred
over $BB$-matches, there is a larger reward for the polymer to $BA$-localize 
(stay close to the interface when diagonally crossing a $B$-block) than to 
$AB$-localize (stay close to the interface when diagonally crossing an $A$-block).

The following is the analogue of Theorem \ref{phtrcurve}, and will be proved in
Section \ref{S4.2.2}. Two constants $0<\alpha_0<\alpha_1<\infty$ appear, which 
will be identified in Section \ref{S2.2}.

\bt{phtrcurve*}
Let $p<p_c$.\\
(i) $\partial\cD$ lies on or below the supercritical curve.\\
(ii) $\partial\cD$ is continuous and intersects each line from the origin 
with slope in $[-1,1)$ at most once.\\ 
(iii) $\partial\cD$ contains the diagonal segment $\{(\alpha,\alpha)\colon\,
\alpha\in [0,\alpha^*]\}$, with $\alpha^*$ the same constant as in Theorem 
{\rm \ref{phtrcurve}(iii)}, but lies below the diagonal elsewhere.\\
(iv) There exists an $\alpha^*(p)\in (0,\infty)$ such that the intersection of
$\partial\cD$ with the lower half of $\CONE$ is the linear segment 
$\{(\beta+\alpha^*(p),\beta)\colon\,\beta\in [-\frac12 \alpha^*(p),0]\}$.\\
(v) As $p \da 0$, $\partial\cD$ converges to the union of the diagonal segment 
in (iii) and the mirror image of the analytic continuation of the supercritical 
curve outside $\CONE$ (i.e., the mirror image of the curved dotted line in Fig. 5). 
In particular, $\lim_{p\da 0} \alpha^*(p)=\alpha_0$.\\
(vi) As $p \ua p_c$, $\partial\cD$ does not converge to the supercritical
curve in Fig.\ 5. In particular, $\lim_{p\ua p_c} \alpha^*(p)=\alpha_1$.
\et


\setlength{\unitlength}{0.5cm}

\begin{picture}(12,12)(-6,-2.5)

  \put(0,0){\line(12,0){12}}
  \put(0,0){\line(0,8){8}}
  \put(0,0){\line(0,-3){3}}
  \put(0,0){\line(-3,0){3}}
  {\thicklines
   \qbezier(0,0)(2,2)(4,4)
   \qbezier(4,4)(3.5,2)(2.5,0)
   \qbezier(2.5,0)(1.875,-.625)(1.25,-1.25)
  }
  \qbezier[40](.75,-.7)(3.5,1.5)(4,4)
  \qbezier[20](4,0)(4,2)(4,4)
  \qbezier[60](4,4)(6,6)(8,8) 
  \qbezier[60](4,4)(4.5,5)(5,8)
  \qbezier[40](-3,3)(0,0)(3.2,-3.0)
  \put(-.8,-.8){$0$}  
  \put(12.5,-0.2){$\alpha$}
  \put(-0.1,8.5){$\beta$}
  \put(4.3,.3){$\alpha^*$}
  \put(2.3,-.9){$\alpha^*(p)$}
  \put(4,4){\circle*{.25}}
  \put(4,0){\circle*{.25}}
  \put(2.55,0){\circle*{.25}}
  \put(6.5,1.4){$\cal L$}
  \put(6.5,-1.7){$\cal L$}
  \put(1.7,.9){$\cal D$}

  \put(-4,-5){\small 
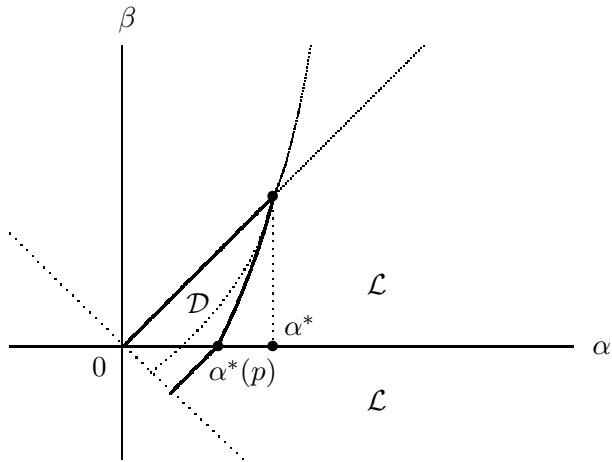
               Fig.\ 7. Qualitative picture of $\partial\cD$ for $p < p_c$. 
               The curved dotted line
               \normalsize}
  \put(-4,-5.8){\small
              is the mirror image of the union of the supercritical curve
              off the
              \normalsize}
  \put(-4,-6.6){\small
              diagonal and its analytic continuation outside $\CONE$
              (see Fig.\ 5).
              \normalsize} 

\end{picture}

\vskip 2.5truecm


It is clear from (\ref{DL*}) that the part off the diagonal is a critical line. 
We will see in Section \ref{S4.2.3} that also the part on the diagonal is a 
critical line.

We will see in Section \ref{S2.2} that $\cD$ corresponds to the situation where
the polymer is \emph{fully delocalized into the $A$-blocks and the $B$-blocks}, 
while $\cL$ corresponds to the situation where the polymer is \emph{partially 
$BA$-localized}. We expect that $\cD$ is strictly increasing in $p$ and that 
$\alpha^*>\alpha_1$, but we are unable to prove this. The curved dotted line
crosses the horizontal axis at $\alpha_0$. 

We will see in Section \ref{S4.3} that $\cL$ contains a second curve (see Fig.\ 8)
at which a phase transition occurs from \emph{partially $BA$-localized} to 
\emph{partially $BA$-localized and partially $AB$-localized}. Qualitatively, 
this curve behaves like the supercritical curve (e.g.\ it also starts at the 
point $(\alpha^*,\alpha^*)$), but unfortunately we know little about it. We 
expect it to be strictly increasing in $\alpha$. We expect it to move down 
as $p$ increases. We do know that it converges to the supercritical curve as 
$p\ua p_c$.


\setlength{\unitlength}{0.5cm}

\begin{picture}(12,12)(-6,-3)

  \put(0,0){\line(10,0){10}}
  \put(0,0){\line(0,7){7}}
  \put(0,0){\line(0,-2){2}}
  \put(0,0){\line(-2,0){2}}
  {\thicklines
   \qbezier(4,4)(6,5)(10,6)
  }
  \qbezier[60](0,0)(3.5,3.5)(7,7)
  \qbezier[30](4,0)(4,2)(4,4)
  \qbezier[60](0,6.5)(5,6.5)(10,6.5)
  \qbezier[40](4,4)(3.5,1.5)(2.5,0)
  \qbezier[10](2.5,0)(2,-.5)(1.35,-1.15)
  \qbezier[40](-2,2)(0,0)(2,-2)
  \put(-.8,-.8){$0$}  
  \put(10.5,-0.2){$\alpha$}
  \put(-0.1,7.5){$\beta$}
  \put(3.8,-.8){$\alpha^*$}
  \put(4,4){\circle*{.25}}
  \put(6.5,4.2){$\cal L$}
  \put(6.5,5.5){$\cal L$}
  \put(2,1){$\cal D$}
  
  \put(-4,-3.5){\small 
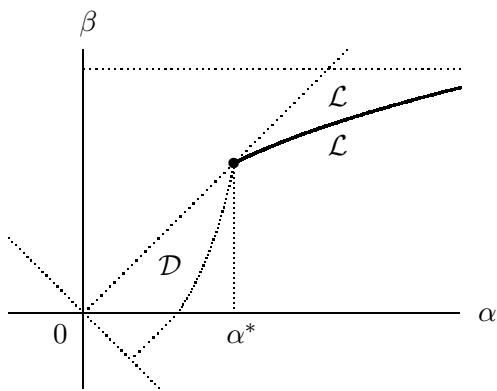
               Fig.\ 8. Conjectured critical line inside $\cL$ for $p < p_c$. 
               \normalsize}
\end{picture}

\vskip 1truecm


To prove Theorem \ref{phtrcurve*}, we will \emph{reformulate} the criteria 
$\psi_{AB}(\bar x)>\psi_{AA}(\bar x)$ and $\psi_{BA}(\bar y)>\psi_{BB}(\bar y)$  
in terms of criteria for the free energy of a model with a \emph{single linear 
interface}. This reformulation, which will be given in Section \ref{S2.4}, is 
again crucial in allowing us to get a handle on the critical curve in Fig.\ 7.

\subsection{Heuristic explanation of the phase diagram}
\label{S1.5extra}

The physical background of the three critical curves in Figs.\ 5, 7 and 8
is as follows. 

\medskip\noindent
$\bullet$ $p \geq p_c$:

\medskip
Consider the boundary $\partial\cD$ sketched in Fig.\ 5. Pick a point 
$(\hat\alpha,\hat\beta$) inside $\cD$. Then, since $p \geq p_c$ and 
$\alpha\geq\beta$, the polymer spends almost all of its time deep 
inside $A$-blocks. Now increase $\beta$ but keep $\alpha=\hat\alpha$ 
fixed. Then there will be a larger energetic advantage for the polymer 
to move some of its monomers from the $A$-blocks to the $B$-blocks by 
crossing the interface inside the $AB$-block pairs. There is some 
entropy loss associated with doing so. The polymer has three options: 
(i) it may place all $A$-monomers in the $A$-blocks and no monomers 
in the $B$-blocks (resulting in all energy coming from $AA$-matches 
and some entropy); (ii) it may place all $A$-monomers in the $A$-blocks 
and a positive fraction of $B$-monomers in the $B$-blocks (resulting 
in a higher energy and a lower entropy); (iii) it may sacrifice some 
fraction of $AA$-matches to get a larger fraction of $BB$-matches
(resulting in an even higher energy and an even lower entropy). If 
$\beta$ is large enough, then the energy advantage will dominate, so 
that $AB$-localization sets in. The value at which this happens depends 
on $\hat\alpha$ and is strictly positive. Since the entropy loss is finite, 
for $\hat\alpha$ large enough the energy-entropy competition plays out 
not only below the diagonal, but also below a horizontal asymptote. The 
larger the value of $\hat\alpha$, the larger the value of $\beta$ where 
$AB$-localization sets in. This explains why the part of $\partial\cD$
off the diagonal moves to the right and up.

\medskip\noindent
$\bullet$ $p < p_c$:

\medskip
First consider the boundary $\partial\cD$ sketched in Fig.\ 7. Pick a point 
$(\hat\alpha,\hat\beta$) inside $\cD$. Since $p<p_c$, the polymer spends
almost all of its time deep inside $A$-blocks and $B$-blocks. Now increase 
$\alpha$ but keep $\beta=\hat\beta$ fixed. Then, while remaining delocalized, 
the polymer will spend more time in the $A$-blocks and less time in the 
$B$-blocks, trying to lower its energy with some attendant loss of entropy. 
As $\alpha$ increases further, there will be a larger energetic advantage 
for the polymer to move some of its monomers from the $B$-blocks to the 
$A$-blocks by crossing the interface inside the $BA$-block pairs. If $\alpha$ 
is large enough, then the energetic advantage will dominate, so that 
$BA$-localization sets in eventually. The value of $\alpha$ at which this 
happens depends on $\hat\beta$. A larger value of $\hat\beta$ means that 
the polymer spends more time in the $B$-blocks (at fixed $\alpha$) with 
larger entropy. Consequently, more entropy will be lost on $BA$-localization 
and the value of $\alpha$ where $BA$-localization sets in will be larger. This
explains why the part of $\partial\cD$ off the diagonal moves to the right 
and up. Similarly, if $p$ decreases, then the polymer hits more $B$-blocks, 
and to compensate for the loss of energy it will spend more time in an $A$-block 
when it hits one and less time in a $B$-block when it hits one (at fixed $\hat\alpha$ 
and $\hat\beta$). Consequently, less entropy will be lost on $BA$-localization and 
the value of $\alpha$ (at fixed $\hat\beta$) where $BA$-localization sets in will 
be smaller. This explains why $\cD$ shrinks with $p$.

If $\beta\leq 0$, then there is a penalty for having $B$-monomers in 
$B$-blocks. Therefore, when the polymer $BA$-localizes, it will spend 
all the time it runs along the interface in the $A$-block and then shoot 
through the interface to spend its remaining time in the $B$-block (on 
its way to the diagonally opposite corner). Hence, the energy-entropy 
competition only depends on the difference $\alpha-\beta$. This explains 
why there is a degeneration of the critical curve into a linear segment. 

In the limit as $p \da 0$, the density of $A$-blocks tends to 
zero and so the polymer spends more and more of its time in $B$-blocks. 
Therefore the localization mechanism looks more and more like that for the 
supercritical curve with $\alpha \leftrightarrow \beta$ and $p \ua 1$.

Next consider the curve separating $\cL$ sketched in Fig.\ 8. Pick a point 
$(\hat\alpha,\hat\beta$) inside $\cL$. Now increase $\beta$ but keep 
$\alpha=\hat\alpha$ fixed. Then, as before, an energy-entropy competition 
sets in. The polymer has the same three options inside $AB$-blocks as in the 
supercritical case, and therefore the curve has the same qualitative behaviour. 
In the limit as $p \ua p_c$, the polymer spends more and more of its time in 
$A$-blocks. Therefore the $AB$-localization mechanism looks more and more like 
that for the supercritical curve. If $p$ decreases, then the polymer hits more 
$B$-blocks, and to compensate for the loss of energy it will spend more time 
in an $A$-block when it hits one and less time in a $B$-block when it hits one 
(at fixed $\hat\alpha$ and $\hat\beta$). Consequently, more entropy will be 
lost on $AB$-localization and the value of $\beta$ (at fixed $\hat\alpha$) where 
$AB$-localization sets in will be larger. This explains why the curve moves up 
as $p$ decreases.

\bigskip
Finally, with the help of the two symmetry properties stated in (\ref{symms}), 
the phase diagram can be extended from $\CONE$ to $\R^2$. When doing so, we
obtain the following phase diagram. Here, the label on $\cD$ ($\cL$) indicates 
the type of (de)localization.


\vskip 1truecm

\setlength{\unitlength}{0.65cm}

\begin{picture}(10,10)(-8,-3)

  \qbezier[60](-7,0)(0,0)(7,0)
  \qbezier[60](0,-7)(0,0)(0,7)
  {\thicklines
  \qbezier(-2,-2)(0,0)(2,2)
  \qbezier(2,2)(4,4)(6,6)
  \qbezier(2,2)(4,3)(7,4)
  \qbezier(2,2)(3,4)(4,7)
  \qbezier(-2,-2)(-4,-4)(-6,-6)
  \qbezier(-2,-2)(-4,-3)(-7,-4)
  \qbezier(-2,-2)(-3,-4)(-4,-7)
  \qbezier(2,2)(1,2)(0,1.7)
  \qbezier(-2,-2)(-2,-1)(-1.7,0)
  \qbezier(0,1.7)(-.85,.85)(-1.7,0) 
  }
  \put(7.5,-0.2){$\alpha$}
  \put(-0.1,7.5){$\beta$}
  \put(2,2){\circle*{.25}}
  \put(-2,-2){\circle*{.25}}
  \put(0,1.7){\circle*{.25}}
  \put(-1.7,0){\circle*{.25}}

  \put(3,.7){\tiny $\cD_A$}
  \put(3,-.7){\tiny $\cD_A$}
  \put(.5,-3){\tiny $\cD_A$}
  \put(-1.1,-3){\tiny $\cD_A$}

  \put(-.8,.5){\tiny $\cD_{A+B}$}

  \put(4.8,4){\tiny $\cL_{AB}$}
  \put(-4.5,-5){\tiny $\cL_{AB}$}

  \put(.5,3){\tiny $\cL_{AB}$}
  \put(-1.4,3){\tiny $\cL_{AB}$}
    
  \put(-4,.7){\tiny $\cL_{AB}$}
  \put(-4,-.7){\tiny $\cL_{AB}$}

  \put(5,7){\tiny $\cL_{AB+BA}$}
  \put(-8.5,-5.2){\tiny $\cL_{AB+BA}$}
    
  \put(-2,-8.5){\small 
               Fig.\ 9. Full phase diagram for $p\geq p_c$. 
               \normalsize}
\end{picture}

\vskip 4.5truecm



\setlength{\unitlength}{0.65cm}

\begin{picture}(10,10)(-8,-3)

  \qbezier[60](-7,0)(0,0)(7,0)
  \qbezier[60](0,-7)(0,0)(0,7)
  {\thicklines
  \qbezier(-2,-2)(0,0)(2,2)
  \qbezier(2,2)(4,3)(7,4)
  \qbezier(2,2)(3,4)(4,7)
  \qbezier(-2,-2)(-4,-3)(-7,-4)
  \qbezier(-2,-2)(-3,-4)(-4,-7)
  \qbezier(2,2)(1,2)(0,1.7)
  \qbezier(-2,-2)(-2,-1)(-1.7,0)
  \qbezier(0,1.7)(-.85,.85)(-1.7,0)
  
  \qbezier(2,2)(2,1)(1.7,0)
  \qbezier(-2,-2)(-1,-2)(0,-1.7)
  \qbezier(1.7,0)(.85,-.85)(0,-1.7)

  }
  \put(7.5,-0.2){$\alpha$}
  \put(-0.1,7.5){$\beta$}
  \put(2,2){\circle*{.25}}
  \put(-2,-2){\circle*{.25}}
  \put(0,1.7){\circle*{.25}}
  \put(-1.7,0){\circle*{.25}}

  \put(1.7,0){\circle*{.25}}
  \put(0,-1.7){\circle*{.25}}

  \put(3,.6){\tiny $\cL_{BA}$}
  \put(3,-.7){\tiny $\cL_{BA}$}
  \put(.5,-3){\tiny $\cL_{BA}$}
  \put(-1.3,-3){\tiny $\cL_{BA}$}

  \put(-.8,.5){\tiny $\cD_{A+B}$}
  \put(-.3,-.7){\tiny $\cD_{A+B}$}

   \put(.5,3){\tiny $\cL_{AB}$}
  \put(-1.4,3){\tiny $\cL_{AB}$}
    
  \put(-4,.6){\tiny $\cL_{AB}$}
  \put(-4,-.7){\tiny $\cL_{AB}$}

  \put(4.5,5){\tiny $\cL_{AB+BA}$}
  \put(-6.5,-5.2){\tiny $\cL_{AB+BA}$}
  
  \put(-2,-8.5){\small 
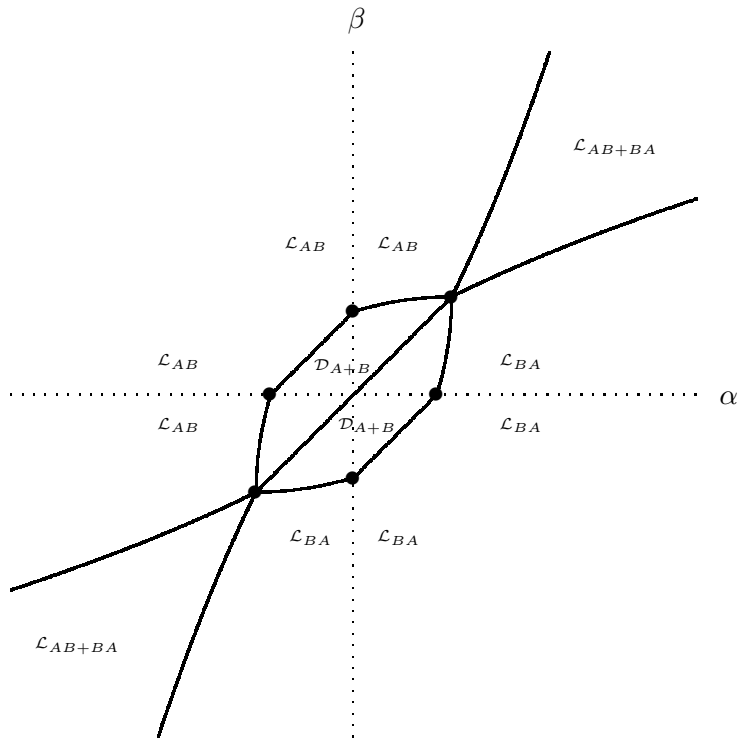
               Fig.\ 10. Full phase diagram for $1-p_c<p<p_c$. 
               \normalsize}
\end{picture}

\vskip 4truecm


\noindent
The figure for $p \leq 1-p_c$ is the same as for $p \geq p_c$, but with all 
the phases reflected in the first diagonal and with the labels $A$ and $B$ 
interchanged. Note that the phase diagram is discontinuous both at $p=p_c$ 
and $p=1-p_c$.

\subsection{Open problems}
\label{S1.6}

The fine details of the two subcritical curves remain to be settled.
Here are some further open problems:

\begin{itemize}
\item[1.] 
Are the critical curves smooth off the diagonal and inside the first quadrant? 
Even for the model with a single linear interface this question has not been 
settled. 
\item[2.]
For $p \geq p_c$, is the free energy infinitely differentiable inside
the localized phase? For the model with a single linear interface this
was proved by Giacomin and Toninelli \cite{GT05b}. Is the same true
for $p<p_c$ in the interior of the two subphases of the localized phase?
\item[3.]
Our phase transitions are defined in terms of a non-analyticity in the 
free energy. Heuristically, they correspond to the path changing its 
behaviour from being fully delocalized away from the interfaces to being 
partially localized near the interfaces (in the subcritical case even in
two possible ways). How can we prove that the path actually has this 
behaviour under the transformed path measure
\be{pathmeas}
P^{\omega,\Omega}_{n,L_n}(\pi) 
= \frac{1}{Z^{\omega,\Omega}_{n,L_n}}\,\,
\exp\left[-H^{\omega,\Omega}_{n,L_n}(\pi)\right]
\ee
for large $n$? For the model with a single linear interface this question 
was settled in Biskup and den Hollander \cite{BdH99} and in Giacomin and
Toninelli \cite{GT05}.
\item[4.]
How does the free energy behave near the critical curve? For the 
model with a single linear interface it was shown by Giacomin and Toninelli
\cite{GT05a} that the phase transition is at least of second order. Numerical 
results in Causo and Whittington \cite{CW03} suggest that the same is true for 
the self-avoiding walk model.
\item[5.]
The coarse-graining expressed by (\ref{Ln}) and the restriction that the polymer 
can enter and exit a pair of neighbouring blocks only at diagonally opposite 
corners are necessary to make the model mathematically tractable. Indeed, the 
corner restriction and $L_n\to\infty$ guarantee that the polymer ``sees one pair 
of blocks at a time'' and self-averages in $\omega$ in each block, which is why 
the free energy can be decomposed into contributions coming from single pairs of 
blocks, while $L_n/n\to 0$ guarantees that the polymer ``sees many blocks'' and 
self-averages in $\Omega$, which is why percolation effects enter. What happens 
when we remove the corner restriction? What happens when the blocks have random 
sizes? 
\end{itemize}


\section{Preparations}
\label{S2}

In Section \ref{S2.1} we compute entropies for paths that cross a block and paths
that run along an interface. In Section \ref{S2.2} we derive a formula for $\psi_{kl}$
in (\ref{fekl}) and $S_{kl}$ in (\ref{Sdef}). In Section \ref{S2.3} we deduce a
criterion for localization when $p \geq p_c$ \emph{in terms of the free energy
for the model with a single linear interface}. In Section \ref{S2.4} we do the 
same when $p<p_c$. In Section \ref{S2.5} we analyse the supremum $\rho^*(p)$ in 
(\ref{rho*def}) and the variational formula for $F(\rho)$ in (\ref{Fdef}).

\subsection{Path entropies}
\label{S2.1}

The results in this section are based on straightforward computations, but are
crucial for the rest of the paper.

\subsubsection{Paths crossing a block}
\label{S2.1.1}

An important ingredient in the identification of $\psi_{kl}(a)$, $k,l\in\{A,B\}$, 
is the following combinatorial lemma. Let 
\be{DOMdef}
\DOM = \{(a,b)\colon\, a \geq 1+b, b\geq 0\}.
\ee 
For $(a,b) \in \DOM$, let $N_L(a,b)$ denote the number of $aL$-step self-avoiding 
directed paths from $(0,0)$ to $(bL,L)$ whose vertical displacement stays within 
$(-L,L]$ ($aL$ and $bL$ are integer). Let
\be{kappa}
\kappa(a,b) = \lim_{L\to\infty} \frac{1}{aL} \log N_L(a,b).
\ee

\bl{l:ka}
(i) $\kappa(a,b)$ exists and is finite for all $(a,b)\in\DOM$.\\
(ii) $(a,b) \mapsto a\kappa(a,b)$ is continuous and strictly concave
on $\DOM$ and analytic on the interior of $\DOM$.\\ 
(iii) For all $a \geq 2$,
\be{ka}
a \kappa(a,1) = \log 2 + \frac 12 \left[a \log a - (a-2) \log (a-2)\right]. 
\ee
(iv) $\sup_{a\geq 2} \kappa(a,1) = \kappa(a^*,1) = \frac12 \log 5$ with unique maximiser
$a^* = \frac52$.\\
(v) $(\frac{\partial}{\partial a}\kappa)(a^*,1)=0$ and $a^*(\frac{\partial}{\partial b}
\kappa)(a^*,1)=\frac12\log\frac95$.
\el

\bpr
(i) First we do the computation without the restriction on the vertical displacement.
Later we show that putting in the restriction is harmless. 

Let $N^0_L(a,b)$ denote the number of $aL$-step self-avoiding directed paths 
from $(0,0)$ to $(bL,L)$. Since such paths make $bL$ steps to the right, 
$\frac{a+1-b}{2}L$ steps upwards and $\frac{a-1-b}{2}L$ steps downwards, we have
\be{kacomb}
N^0_L(a,b) 
= \sum_{k=1}^{bL} \left(\ba{c} bL\\k \ea\right)
                           \left(\ba{c} \frac{a+1-b}{2}L-1\\k-1 \ea\right)
\sum_{l=1}^{bL-k} \left(\ba{c} bL-k\\l \ea\right)
                             \left(\ba{c} \frac{a-1-b}{2}L-1\\l-1 \ea\right).
\ee
Here, $k$ counts the number of columns where the path moves upward, $l$ counts the
number of colums where the path moves downward, the first and the third binomial 
coefficient count the number of choices for these columns, while the second and
the fourth binomial coefficient count the number of ways in which the prescribed
number of steps can be distributed over these columns. Since, by Stirling's formula,
\be{Stir}
\lim_{L\to\infty} \frac{1}{L} \log \left(\ba{l} uL\\vL \ea\right)
= u \log u - v \log v - (u-v) \log (u-v), \qquad 0 \leq v \leq u,
\ee
we get, by putting $k=\delta L$, $l=\epsilon L$, that
\be{max}
a \kappa^0 (a,b) = \lim_{L\to\infty} \frac{1}{L} \log N^0_L(a,b)
= \sup_{\delta,\epsilon} f_{ab}(\delta,\epsilon)
\ee
with
\be{fdef}
\begin{aligned}
f_{ab}(\delta,\epsilon) &= b \log b - 2 \delta \log \delta
+ \left(\frac{a+1-b}{2}\right) \log \left(\frac{a+1-b}{2}\right)\\
&\qquad - \left(\frac{a+1-b}{2}-\delta\right) \log \left(\frac{a+1-b}{2}-\delta\right)
- 2 \epsilon \log \epsilon -(b-\delta-\epsilon) \log (b-\delta-\epsilon)\\
&\qquad  + \left(\frac{a-1-b}{2}\right) \log \left(\frac{a-1-b}{2}\right)
- \left(\frac{a-1-b}{2}-\epsilon\right) \log \left(\frac{a-1-b}{2}-\epsilon\right).
\end{aligned}
\ee
Computing
\be{fdif}
\begin{aligned}
\frac{\partial f_{ab}}{\partial\delta} 
&= \log \left[\frac{(\frac{a+1-b}{2}-\delta)(b-\delta-\epsilon)}{\delta^2}\right],\\
\frac{\partial f_{ab}}{\partial\epsilon} 
&= \log \left[\frac{(\frac{a-1-b}{2}-\epsilon)(b-\delta-\epsilon)}{\epsilon^2}\right],  
\end{aligned}
\ee
and setting these derivatives equal to zero, we find that the maximisers $\delta_{ab}$ 
and $\epsilon_{ab}$ of the right-hand side of (\ref{max}) are solutions of quadratic 
equations, namely, 
\be{quadr}
\begin{aligned}
0 &= (1+b)\delta^2 - (a+1)b\delta +\frac{a+1-b}{2}\,b^2,\\
0 &= (1-b)\epsilon^2 +(a-1)b\epsilon -\frac{a-1-b}{2}\,b^2,
\end{aligned}
\ee
which leads to 
\be{maximisers}
\begin{aligned}
\delta_{ab} &= \frac{b}{2(1+b)}\left[(a+1)-\sqrt{(a-b)^2+(b^2-1)}\right],\\
\epsilon_{ab} &= \frac{b}{2(1-b)}\left[-(a-1)+ \sqrt{(a-b)^2+(b^2-1)}\right],
\end{aligned}
\ee
for $b \neq 1$, and
\be{maximisers*}
\delta_{a1} = \frac12, \qquad \epsilon_{a1} = \frac{a-2}{2(a-1)},
\ee
for $b=1$. Substitution of (\ref{maximisers}--\ref{maximisers*}) into 
(\ref{max}--\ref{fdef}) yields a formula for $a\kappa^0(a,b)$ in closed 
form. From this formula it is obvious that $(a,b)\mapsto\kappa^0(a,b)$ 
is continuous on $\DOM$ and analytic on the interior of $\DOM$.

It remains to show that the restriction on the vertical displacement has no effect 
in the limit as $L\to\infty$. This can be done by appealing to the reflection 
principle. Indeed, let $N^\da_L(a,b)$ be the number of paths where the restriction 
of not moving above the line of height $L$ is inserted. Then $N^\da_L(a,b)$ 
is the difference of two terms of the type $N^0_L(a,b)$ in (\ref{kacomb}), 
one with the path ending at $(bL,L+2)$ and one with the path ending at $(bL,L)$. 
A little computation shows that this difference equals $N^0_L(a,b)$ divided by 
a term that is growing at most polynomially fast in $L$. This polynomial factor does 
not affect the exponential asymptotics. A similar argument shows that the restriction 
of not moving below the line of height $-L+1$ is harmless as well. Hence
\be{kappaeq}
\kappa(a,b)=\kappa^0(a,b).
\ee

\noindent
(ii) Pick any $b_1,b_2 \geq 0$ and $a_1\geq 1+b_1,a_2\geq 1+b_2$. Consider a block
of height $L$ and width $\frac12(b_1+b_2)L$, and partition this block into four parts
by cutting it at height $\frac12 L$ and width $\frac12 b_1L$. The number of paths
that cross the large block in $\frac12(a_1+a_2)L$ steps is larger than or equal
to the number of paths that cross the lower left block in $\frac12 a_1L$ steps
times the number of paths that cross the upper right block in $\frac12 a_2L$ steps,
i.e.,
\be{convN}
N^0_L\left(\frac{a_1+a_2}{2},\frac{b_1+b_2}{2}\right) 
\geq N^0_{\frac12 L}(a_1,b_1) N^0_{\frac12 L}(a_2,b_2). 
\ee
By (\ref{max}) and (\ref{kappaeq}), this proves that
\be{kappaconv}
\frac{a_1+a_2}{2}\kappa\left(\frac{a_1+a_2}{2},\frac{b_1+b_2}{2}\right)
\geq \frac12 a_1\kappa(a_1,b_1) + \frac12 a_2\kappa(a_2,b_2),
\ee
which is the concavity desired. Strict concavity follows from analyticity on
the interior of $\DOM$, because $a\kappa(a,b)$ clearly is not linear in either
$a$ or $b$.

\noindent 
(iii) Substitute (\ref{maximisers*}) into (\ref{max}--\ref{fdef}) to get 
the formula for $a\kappa(a,1)$ stated in (\ref{ka}).

\noindent
(iv) Since $\frac{d}{da}\kappa(a,1)=-\frac{1}{a^2}\log[2(a-2)]$, the supremum
is uniquely attained at $a^*=\frac52$, giving the claim.

\noindent
(v) Compute, from (\ref{fdef}),
\be{parak}
\begin{aligned}
\left(\frac{\partial}{\partial a}\kappa\right)(a,b)
&= \left(\frac{\partial}{\partial a}\left(\frac{1}{a}f_{ab}\right)\right)
(\delta_{ab},\epsilon_{ab})\\
&= - \frac{1}{a^2}f_{ab}(\delta_{ab},\epsilon_{ab})
+ \frac{1}{a}\left(\frac{\partial}{\partial a}f_{ab}\right)(\delta_{ab},\epsilon_{ab})\\
&= -\frac{1}{a}\kappa(a,b) 
+ \frac{1}{a}\,\frac12\log\left[\frac{(\frac{a+1-b}{2})(\frac{a-1-b}{2})}
{(\frac{a+1-b}{2}-\delta_{ab})(\frac{a-1-b}{2}-\epsilon_{ab})}\right]
\end{aligned}
\ee
and
\be{parbk}
\begin{aligned}
\left(\frac{\partial}{\partial b}\kappa\right)(a,b)
&= \left(\frac{\partial}{\partial b}\left(\frac{1}{a}f_{ab}\right)\right)
(\delta_{ab},\epsilon_{ab})\\
&= \frac{1}{a}\left(\frac{\partial}{\partial b}f_{ab}\right)(\delta_{ab},\epsilon_{ab})\\
&= \frac12\,\frac{1}{a}\log\left[\frac{b^2(\frac{a+1-b}{2}
-\delta_{ab})(\frac{a-1-b}{2}-\epsilon_{ab})}
{(b-\delta_{ab}-\epsilon_{ab})^2(\frac{a+1-b}{2})(\frac{a-1-b}{2})}\right].
\end{aligned}
\ee
Setting $a=a^*=\frac52$, $b=1$, $\delta_{ab}=\delta_{a^*1}=\frac12$ and $\epsilon_{ab}
=\epsilon_{a^*1}=\frac16$, we get the claim. 
\epr 

\subsubsection{Paths running along an interface}
\label{S2.1.2}

We also need the following analogue of Lemma \ref{l:ka}. For $\mu \geq 1$, let 
$\hat N_L(\mu)$ denote the number of $\mu L$-step self-avoiding paths from $(0,0)$ 
to $(L,0)$ with no restriction on the vertical displacement ($\mu L$ is integer). 
Let
\be{kappamudef}
\hat \kappa(\mu) = \lim_{L\to\infty} \frac{1}{\mu L} \log \hat N_L(\mu).
\ee

\bl{l:kamu}
(i) $\hat\kappa(\mu)$ exists and is finite for all $\mu\geq 1$.\\
(ii) $\mu\mapsto\mu\hat\kappa(\mu)$ is continuous and strictly concave
on $[1,\infty)$ and analytic on $(1,\infty)$.\\
(iii) $\hat\kappa(1)=0$ and $\mu\hat\kappa(\mu)\sim\log\mu$ as $\mu\to\infty$.\\
(iv) $\sup_{\mu\geq 1} \mu[\hat\kappa(\mu)-\frac12\log 5]<\frac12\log\frac95$.
\el

\bpr
(i) Similarly as in (\ref{kacomb}),
\be{kacombext}
\hat N_L(\mu) 
= \sum_{k=1}^L \left(\ba{c} L\\k \ea\right)
                           \left(\ba{c} \frac{\mu-1}{2}L-1\\k-1 \ea\right)
\sum_{l=1}^{L-k} \left(\ba{c} L-k\\l \ea\right)
                             \left(\ba{c} \frac{\mu-1}{2}L-1\\l-1 \ea\right).
\ee
Again putting $k=\delta L$, $l=\epsilon L$, we get
\be{maxext}
\mu\hat\kappa(\mu) = \lim_{L\to\infty} \frac{1}{L} \log \hat N_L(\mu)
= \sup_{\delta,\epsilon} f_\mu(\delta,\epsilon)
\ee
with
\be{fdefext} 
\begin{aligned}
f_\mu(\delta,\epsilon) &= -2\delta\log\delta - 2\epsilon\log\epsilon
- (1-\delta-\epsilon)\log(1-\delta-\epsilon)\\
&\qquad -\left(\frac{\mu-1}{2}-\delta\right)\log\left(\frac{\mu-1}{2}-\delta\right)
-\left(\frac{\mu-1}{2}-\epsilon\right)\log\left(\frac{\mu-1}{2}-\epsilon\right)\\
&\qquad + (\mu-1)\log\left(\frac{\mu-1}{2}\right).
\end{aligned}
\ee
Computing
\be{fdifext}
\begin{aligned}
\frac{\partial f_\mu}{\partial\delta} 
&= \log \left[\frac{(\frac{\mu-1}{2}-\delta)(1-\delta-\epsilon)}{\delta^2}\right],\\
\frac{\partial f_\mu}{\partial\epsilon} 
&= \log \left[\frac{(\frac{\mu-1}{2}-\epsilon)(1-\delta-\epsilon)}{\epsilon^2}\right],  
\end{aligned}
\ee
and setting these derivatives equal to zero, we find that the maximisers 
$\delta_\mu$ and $\epsilon_\mu$ of the right-hand side of (\ref{maxext}) are 
equal, $\delta_\mu=\epsilon_\mu$, with $\delta_\mu$ the solution of the quadratic 
equation
\be{quadrext}
0 = \delta^2 -\mu\delta + \frac{\mu-1}{2},
\ee
which leads to
\be{maximisersext}
\delta_\mu = \frac12 \left[\mu - \sqrt{(\mu-1)^2+1}\right].
\ee 
Substitution of (\ref{maximisersext}) into (\ref{fdefext}) yields a formula for 
$\mu\hat\kappa(\mu)$ in closed form. From this formula it is obvious that
$\mu\mapsto\hat\kappa(\mu)$ is continuous on $[1,\infty)$ and analytic on
$(1,\infty)$. 

\noindent
(ii) Pick any $\mu_1,\mu_2\geq 1$. The number of $\frac12(\mu_1+\mu_2)L$-step
paths from $(0,0)$ to $(L,0)$ is larger than or equal to the number of $\frac12 
\mu_1L$-step paths from $(0,0)$ to $(\frac12 L,0)$ times the number of $\frac12 
\mu_2L$-step paths from $(\frac12 L,0)$ to $(L,0)$, i.e.,
\be{convN*}
\hat N_L\left(\frac{\mu_1+\mu_2}{2}\right) 
\geq \hat N_{\frac12 L}(\mu_1)\hat N_{\frac12 L}(\mu_2).
\ee
Via (\ref{kappamudef}), this proves that
\be{kappaconv*}
\frac{\mu_1+\mu_2}{2}\hat\kappa\left(\frac{\mu_1+\mu_2}{2}\right)
\geq \frac12\mu_1\hat\kappa(\mu_1) + \frac12\mu_2\hat\kappa(\mu_2),
\ee
which is the concavity desired. Strict concavity follows from smoothness on 
$(1,\infty)$, because $\mu\kappa(\mu)$ clearly is not linear in $\mu$.

\noindent
(iii) From (\ref{maximisersext}) we see that $\delta_1(=\epsilon_1)=0$. Hence 
(\ref{maxext}--\ref{fdefext}) give $\hat\kappa(1)=0$. Similarly, if $\mu\to\infty$,
then $\delta_\mu=\frac12 [1-\frac{1}{2\mu}+O(\frac{1}{\mu^2})]$ and hence 
$\mu\hat\kappa(\mu)\sim\log\mu$.

\noindent
(iv) For any $a \geq 2$, $0<b\leq 1$, $\mu\geq 1$ such that $(\mu-1)b \leq a-2$, 
we have
\be{kapparels}
a\kappa(a,1) \geq b\mu\hat\kappa(\mu) + (a-b\mu)\kappa(a-b\mu,1-b).
\ee
Indeed, any $aL$-step self-avoiding path from $(0,0)$ to $(L,L)$ may follow
the interface over a distance $bL$ during $b\mu L$ steps and then wander away
from the interface to the diagonally opposite corner over a distance $(1-b)L$
during $(a-b\mu)L$ steps (see Fig.\ 4). Rewrite (\ref{kapparels}) as
\be{kapparels1} 
\mu\hat\kappa(\mu) \leq 
\frac{1}{b}\left[a\kappa(a,1)-(a-b\mu)\kappa(a-b\mu,1-b)\right].
\ee
Pick $a=a^*$ and let $b\da 0$, to obtain
\be{kapparels2}
\mu\hat\kappa(\mu) \leq \mu\left(\frac{\partial}{\partial a}(a\kappa)\right)(a^*,1)
+\left(\frac{\partial}{\partial b}(a\kappa)\right)(a^*,1).
\ee
By Lemma \ref{l:ka}(iv,v), the right-hand side equals $\mu\frac12\log 5+\frac12
\log\frac95$. Since $\mu\geq 1$ is arbitrary, this proves that $\sup_{\mu\geq 1}
\mu[\hat\kappa(\mu)-\frac12\log 5] \leq \frac12\log\frac95$, which is the claim 
with $\leq$ instead of $<$. A calculation with MAPLE gives that the supremum in 
the left-hand side is attained at $\mu\approx 2.12$ and equals $\approx 0.16$. 
The right-hand side equals $0.29$.
\epr

In Section \ref{S4.1} we will need two special values of $\alpha$, namely,
$\alpha_0$ and $\alpha_1$ given by
\be{alpha01}
\begin{aligned}
\sup_{\mu\geq 1} \mu\left[\hat\kappa(\mu)+\frac12\alpha_0-\frac12\log 5\right]
&= \frac12\log\frac95,\\
\sup_{\mu\geq 1} \mu\left[\hat\kappa(\mu)-\frac12\log 5\right]
&= \frac12\log\left[\frac{4e^{-\alpha_1}(5+e^{-\alpha_1})^2}
{5(5-e^{-\alpha_1})^2}\right].
\end{aligned}
\ee
It follows from Lemma \ref{l:kamu}(iii-iv) that $\alpha_0,\alpha_1>0$. A 
calculation with MAPLE gives the values
\be{alphavals}
\alpha_0 \approx 0.125, \qquad \alpha_1 \approx 0.154.
\ee

\subsection{Free energies per pair of blocks}
\label{S2.2}

In this section we identify $S_{kl}=S_{kl}(\alpha,\beta)$.

\subsubsection{Identification of $S_{AA}$ and $S_{BB}$}
\label{S2.2.1}

\bp{p:fe}
For all $(\alpha,\beta)\in\R^2$,
\be{fe}
S_{AA}=\sup_{a\geq 2}\psi_{AA}(a) = \frac{1}{2}\alpha + \frac{1}{2}\log 5, \qquad
S_{BB}=\sup_{a\geq 2}\psi_{BB}(a) = \frac{1}{2}\beta + \frac{1}{2}\log 5.
\ee
\ep

\bpr
Recall (\ref{Hamiltonian}) and (\ref{feklaL}--\ref{fekl}). For any $aL$-step path in 
an $AA$-block, about half of the monomers contribute $\alpha$ to the energy, because 
$\sum_{i=1}^{aL} 1\{\omega_i=A\} = \frac12 aL[1+o(1)]$ $\omega$-a.s.\ as $L\to\infty$, 
while the remaining monomers contribute $0$ to the energy. Hence
\be{ka1psi}
\psi_{AA}(a) = \frac{1}{2}\alpha + \kappa(a,1).
\ee
Now use Lemma \ref{l:ka}(iv) to get the claim for $S_{AA}$. The proof for $S_{BB}$ is 
the same.   
\epr

\subsubsection{Identification of $S_{AB}$ and $S_{BA}$}
\label{S2.2.2}

It is harder to obtain information on $S_{AB}=\sup_{a\geq 2}\psi_{AB}(a)$ and $S_{BA}
=\sup_{a\geq 2}\psi_{A}(a)$, because these embody the effect of the presence of the 
$AB$-interface. We first consider the free energy per step when the path moves in the 
vicinity of a \emph{single linear interface} $\cI$ separating a liquid $A$ in 
the upper halfplane from a liquid $B$ in the lower halfplane including the interface 
itself. To that end, for $a \geq b > 0$, let $\cW_{aL,bL}$ denote the set 
of $aL$-step directed self-avoiding paths starting at $(0,0)$ and ending at $(bL,0)$. 
Define 
\be{feinf}
\psi^{\omega,\cI}_L(a,b) = \frac{1}{aL} \log Z^{\omega,\cI}_{aL,bL}
\ee
with
\be{Zinf}
\begin{aligned}
Z^{\omega,\cI}_{aL,bL} 
&= \sum_{\pi\in\cW_{aL,bL}} \exp\left[-H^{\omega,\cI}_{aL}(\pi)\right],\\
H^{\omega,\cI}_{aL}(\pi) 
&= - \sum_{i=1}^{aL}\Big(\alpha 1\{\omega_i=A,\pi_i>0\}
+\beta 1\{\omega_i=B,\pi_i \leq 0\}\Big),
\end{aligned}
\ee
where $\pi_i>0$ means that the $i$-th step lies in the upper halfplane and
$\pi_i \leq 0$ means that the $i$-th step lies in the lower halfplane or
in the interface.

\bl{l:feinflim}
For all $(\alpha,\beta)\in\R^2$ and $a \geq b > 0$,
\be{fesainf}
\lim_{L\to\infty} \psi^{\omega,\cI}_L(a,b) = \psi^{\cI}(a,b)
= \psi^{\cI}(\alpha,\beta;a,b)
\ee
exists $\omega$-a.s.\ and is non-random.
\el

\bpr
Since the polymer starts and ends at the interface, the proof can be done via 
a standard subadditivity argument in which two pieces of the polymer are 
concatenated (see e.g.\ Bolthausen and den Hollander \cite{BdH97} or Orlandini 
\emph{et al} \cite{OTW00}). Indeed, fix $a$ and $b$. Then, for any $L_1$ and $L_2$,
\be{subad}
Z^{\omega,\cI}_{a(L_1+L_2),b(L_1+L_2)}
\geq Z^{\omega,\cI}_{aL_1,bL_1} Z^{\sigma^{aL_1}\omega,\cI}_{aL_2,bL_2},
\ee
where $\sigma$ is the left-shift acting on $\omega$. Define
\be{psidef}
\Psi^{\omega,\cI}_K(a,b) = \log Z^{\omega,\cI}_{K,(b/a)K}.
\ee
Then, for any $K_1(=aL_1)$ and $K_2(=aL_2)$, 
\be{subad+}
\Psi^{\omega,\cI}_{K_1+K_2}(a,b) \geq \Psi^{\omega,\cI}_{K_1}(a,b)
+ \Psi^{\sigma^{K_1}\omega,\cI}_{K_2}(a,b).
\ee
We can now apply Kingman's superadditive ergodic theorem, noting that
$\frac1K \Psi^{\omega,\cI}_K(a,b)$ is bounded from above, to conclude that
\be{King}
\lim_{K\to\infty} \frac{1}{K} \Psi^{\omega,\cI}_K(a,b) = \psi^{\cI}(a,b)
\ee  
exists $\omega$-a.s.\ and is non-random.
\epr

The relation linking $\psi_{AB}(a)$ to $\psi^{\cI}(a,b)$ is the following.

\bl{l:linkinf}
For all $(\alpha,\beta)\in\R^2$ and $a\geq 2$,
\be{psiinflink}
\begin{aligned}
\psi_{AB}(a) &= \psi_{AB}(\alpha,\beta;a)\\
&= \sup_{0 \leq b \leq 1,\,a_1 \geq b,\,a_2 \geq 2-b,\,a_1+a_2=a}
\frac{a_1\psi^{\cI}(a_1,b)+a_2[\frac{1}{2}\alpha+\kappa(a_2,1-b)]}{a_1+a_2}.
\end{aligned}
\ee
\el

\bpr
The idea behind this relation is that the polymer follows the $AB$-interface
over a distance $bL$ during $a_1L$ steps and then wanders away from the 
$AB$-interface to the diagonally opposite corner over a distance $(1-b)L$ 
during $a_2L$ steps. The optimal strategy is obtained by maximising over 
$b$, $a_1$ and $a_2$ (recall Figure 4).

A formal proof goes as follows. Look at the last time $u$ and the last site $(v,0)$ 
on the $AB$-interface before the polymer wanders off. This allows us to write the 
associated partition sum as
\be{Zrel}
Z^\omega_{AB}(aL,L) = \sum_{v=0}^L \sum_{u=v}^{aL-(2L-v)}
Z^{\omega,\cI}_{u,v} Z^{\sigma^v\omega}_{aL-u,L-v},
\ee
where $Z^{\omega,\cI}_{u,v}$ is the partition sum for the single interface model
to go in $u$ steps from $(0,0)$ to $(v,0)$, and $Z^{\sigma^u\omega}_{aL-u,L-v}$ 
is the partition sum to go in $aL-u$ steps from $(v,0)$ to $(L,L)$ without returning 
to the interface. Rewrite (\ref{Zrel}) as
\be{Zrelcont}
Z^\omega_{AB}(aL,L) = [1+o(1)]\, L^2 \int_0^1 db \int_b^{a-(2-b)} da_1 
\,Z^{\omega,\cI}_{a_1L,bL} Z^{\sigma^{bL}\omega}_{(a-a_1)L,(1-b)L}. 
\ee
{}From Lemmas \ref{l:ka} and \ref{l:feinflim} we know that, as $L\to\infty$,
\be{Zasymp}
\begin{aligned}
\frac{1}{L} \log Z^{\omega,\cI}_{a_1L,bL} 
&= [1+o(1)]\, a_1\psi^{\cI}(a_1,b)
\quad \omega-a.s.,\\
\frac{1}{L} \log Z^{\sigma^{bL}\omega}_{(a-a_1)L,(1-b)L} 
&= [1+o(1)]\, (a-a_1)\left(\frac12 \alpha + \kappa(a-a_1,1-b)\right) 
\quad \omega-a.s.
\end{aligned}
\ee
For the latter, note that $\sigma^{bL}\omega$ changes with $L$. However, this
causes no problem, because the distribution of $\omega$ is invariant under shifts
and the shift length $bL$ is independent of $\omega$. Substitution of (\ref{Zasymp}) 
into (\ref{Zrelcont}), and of the resulting expression into (\ref{feklaL}), yields the 
claim after we put $a_2=a-a_1$. Indeed, the right-hand sides of (\ref{Zasymp}) are 
continuous in $b$ and $a_1$.   
\epr  

By obvious scaling, there exists a function $\phi^{\cI}$ such that
\be{feinfscal}
\psi^{\cI}(a,b) = \phi^{\cI}(a/b).
\ee
Therefore Lemma \ref{l:linkinf} yields the following.

\bp{p:SAB}
For all $(\alpha,\beta)\in\R^2$,
\be{SABinf}
S_{AB} = \sup_{a\geq 2} \psi_{AB}(a) 
= \sup_{0 \leq b \leq 1,\,a_1 \geq b,\,a_2 \geq 2-b}
\frac{a_1\phi^{\cI}(a_1/b)+a_2[\frac{1}{2}\alpha+\kappa(a_2,1-b)]}{a_1+a_2},
\ee
\ep

\bpr
Insert (\ref{feinfscal}) into (\ref{psiinflink}) and take the supremum over $a$.
\epr

This completes the identification of $S_{AB}$. The same formula applies for $S_{BA}$
but with $\alpha$ and $\beta$ interchanged, i.e.,
\be{SBArel}
S_{BA}(\alpha,\beta) = S_{AB}(\beta,\alpha).
\ee
Recall from the remark made below (\ref{fekl}) that the first index labels the type 
of the block that is diagonally crossed, while the second index labels the type 
of the block that appears as its neighbour. 

Note that $\phi^\cI$ is symmetric in $\alpha$ and $\beta$. The asymmetry in (\ref{Zinf}), 
coming from the fact that the interface is labelled $B$ while the polymer starts at the 
interface, is not felt in the limit as $L\to\infty$. Further note that
\be{phirelbds}
\begin{aligned}
&\phi^\cI(\alpha,\beta;\mu) \in \left[\frac12\alpha+\hat\kappa(\mu),
\alpha+\hat\kappa(\mu)\right] 
\quad \forall\, \alpha \geq \beta \geq 0,\\
&\phi^\cI(\alpha,\beta;\mu)  = \frac12\alpha+\hat\kappa(\mu) 
\,\,\qquad\qquad\qquad\forall\, \alpha \geq 0 \geq \beta.
\end{aligned} 
\ee

We close with the following facts.

\bl{psiprop}
Let $k,l\in\{A,B\}$.\\
(i) For all $(\alpha,\beta)\in\R^2$, $a\mapsto a\psi_{kl}(\alpha,\beta;a)$ is 
continuous and concave on $[2,\infty)$.\\
(ii) For all $a\in [2,\infty)$, $\alpha\mapsto\psi_{kl}(\alpha,\beta;a)$ and $\beta
\mapsto\psi_{kl}(\alpha,\beta;a)$ are continuous and non-decreasing on $\R$.
\el

\bpr
(i) The claim is trivial for $k=l$, because of the simple form of $\psi_{AA}$ and
$\psi_{BB}$ (recall (\ref{ka}) and (\ref{ka1psi})). The proof for $k \neq l$
runs as follows. Rewrite (\ref{psiinflink}) as
\be{varrel1}
a\psi_{AB}(a) = \sup_{0\leq b\leq 1,a_1\geq b,a_2\geq 2-b,a_1+a_2=a}
\left\{a_1\psi^{\cI}(a_1,b)+a_2\left[\frac12\alpha+\kappa(a_2,1-b)\right]\right\}.
\ee
From this it follows that
\be{varrel2}
\begin{aligned}
&\frac12 a^1\psi_{AB}(a^1)+\frac12 a^2\psi_{AB}(a^2)
= \sup_{0\leq b^1\leq 1,a^1_1\geq b^1,a^1_2\geq 2-b^1,a^1_1+a^1_2=a^1}
\,\,\sup_{0\leq b^2\leq 1,a^2_1\geq b^2,a^2_2\geq 2-b^2,a^2_1+a^2_2=a^2}\\
&\left\{\frac12 a^1_1\psi^{\cI}(a^1_1,b^1) + \frac12 a^2_1\psi^{\cI}(a^2_1,b^2) 
+ \frac12 (a^1_2+a^2_2)\frac12\alpha + \frac12a^1_2\kappa(a^1_2,1-b^1)  
+ \frac12a^2_2\kappa(a^2_2,1-b^2)\right\}.
\end{aligned}
\ee 
A standard concatenation argument gives
\be{varrel3}
\begin{aligned}
\frac12 a^1_1\psi^{\cI}(a^1_1,b^1)
+ \frac12 a^2_1\psi^{\cI}(a^2_1,b^2) 
&= \bar a_1 \psi^{\cI}(\bar a_1,\bar b),\\ 
\frac12 a^1_2\kappa(a^1_2,1-b^1)  
+ \frac12 a^2_2\kappa(a^2_2,1-b^2)
&= \bar a_2 \kappa(\bar a_2,1-\bar b),
\end{aligned}
\ee
where we abbreviate
\be{varrel4}
\bar a_1 = \frac{a^1_1+a^2_1}{2},\quad \bar a_2 = \frac{a^1_2+a^2_2}{2},\quad
\bar b = \frac{b^1+b^2}{2}.
\ee
Since the double supremum in (\ref{varrel2}) is more restrictive than the single 
supremum over $0\leq\bar b\leq 1$, $\bar a_1\geq \bar b$, $\bar a_2\geq 2-\bar b$, 
$\bar a_1+\bar a_2 = \bar a$, with $\bar a = (a^1+a^2)/2$, it follows from 
(\ref{psiinflink}) and (\ref{varrel2}--\ref{varrel3}) that
\be{varrel5}
\frac12 a^1\psi_{AB}(a^1)+\frac12 a^2\psi_{AB}(a^2)
\leq \bar a \psi_{AB}(\bar a).
\ee
A similar argument applies to $\psi_{BA}$, after replacing $\frac12\alpha$ by 
$\frac12\beta$ and noting that $\psi^{\cI}$ is symmetric in $\alpha$ and $\beta$.

\medskip\noindent
(ii) The claim is again trivial for $k=l$. For $k \neq l$, note that $\psi^{\cI}$ 
has the same property, as is evident from (\ref{feinf}--\ref{fesainf}). Hence the
claim follows from Lemma \ref{l:linkinf}.  
\epr

\subsection{Criterion for $S_{AB}>S_{AA}$}
\label{S2.3}

For all $(\alpha,\beta)\in\R^2$, we have
\be{Sineq}
S_{AB} \geq S_{AA}.
\ee
The following gives us a criterion for when strict inequality occurs. In Section
\ref{S4.1.1} this will be proved to be the criterion for localization when $p \geq p_c$.

\bp{p:phtrinfchar}
$S_{AB}>S_{AA}$ if and only if
\be{phinfcr}
\sup_{\mu \geq 1} \mu[\phi^{\cI}(\mu)-S_{AA}] > \frac12\log\frac95.
\ee
\ep

\bpr
{}From Propositions \ref{p:fe} and \ref{p:SAB}, together with the reparametrisation
$\mu=a_1/b$ and $\nu=a_2/b$, it follows that
\be{diffrew}
S_{AB}-S_{AA} = \sup_{\mu\geq 1,\,\nu\geq 1}
\frac{\mu[\phi^{\cI}(\mu)-S_{AA}]-\nu[\frac12 \log 5 - f(\nu)]}{\mu+\nu}
\ee
with
\be{fnudef}
f(\nu) = \sup_{\frac{2}{\nu+1}\leq b\leq 1} \kappa(b\nu,1-b), \qquad \nu \geq 1.
\ee
Abbreviate $g(\nu)=\nu[\frac12\log 5-f(\nu)]$. Below we will show that
\be{fineq}
\begin{aligned}
&{\rm (i)} \,\,\,g(\nu) > \frac12\log\frac95\,\, \mbox{ for all } \nu \geq 1,\\
&{\rm (ii)} \,\,\lim_{\nu\to\infty} g(\nu) = \frac12\log\frac95.
\end{aligned}
\ee
This will imply the claim as follows. If $\mu[\phi^{\cI}(\mu)-S_{AA}] \leq 
\frac12\log\frac95$ for all $\mu$, then by (i) the numerator in (\ref{diffrew}) 
is strictly negative for all $\mu$ and $\nu$, and so by (ii) the supremum is taken at 
$\nu=\infty$, resulting in $S_{AB}-S_{AA}=0$. On the other hand, if 
$\mu[\phi^{\cI}(\mu)-S_{AA}] > \frac12\log\frac95$ for some $\mu$, 
then, for that $\mu$, by (i) and (ii) the numerator is strictly positive 
for $\nu$ large enough, resulting in $S_{AB}-S_{AA}>0$.

To prove (\ref{fineq}), we will need the following inequality. Abbreviate
$\chi(a,b)=a\kappa(a,b)$. Then by Lemma \ref{l:ka}(ii) we have, for all
$(s,t) \neq (u,v)$ in $\DOM$,
\be{chiconv}
\begin{aligned}
\chi(s,t) - \chi(u,v) 
&= \int_0^1 dw\,\frac{\partial}{\partial w}\,\chi(u+w(s-u),v+w(t-v))\\
&> \left[\frac{\partial}{\partial w}\,\chi(u+w(s-u),v+w(t-v))\right]_{w=1}\\
&= (s-u)\left(\frac{\partial}{\partial a}\,\chi\right)(s,t)
+ (t-v)\left(\frac{\partial}{\partial b}\,\chi\right)(s,t).
\end{aligned} 
\ee

\medskip
To prove (\ref{fineq})(i), put $b=a/\nu$ in (\ref{fnudef}) and use Lemma 
\ref{l:ka}(iv,v) to rewrite the statement in (\ref{fineq})(i) as
\be{rewineq}
\kappa\left(a,1-\frac{a}{\nu}\right) < \kappa(a^*,1) - \frac{a^*}{\nu}
\left(\frac{\partial}{\partial b}\kappa\right)(a^*,1)\,\,
\mbox{ for all } \nu \geq 1 \mbox{ and } \frac{2\nu}{\nu+1}\leq a\leq \nu.
\ee
But this inequality follows from (\ref{chiconv}) by picking $s=a^*$, $t=1$, 
$u=a$, $v=1-\frac{a}{\nu}$, cancelling a term $a^*\kappa(a^*,1)$ on both sides, 
using that $(\frac{\partial}{\partial a}\kappa)(a^*,1)=0$, and afterwards 
cancelling a common factor $a$ on both sides.

\medskip
To prove (\ref{fineq})(ii), we argue as follows. Picking $b=\frac{a^*}{\nu}$ 
in (\ref{fnudef}), we get from Lemma \ref{l:ka}(iv) that 
\be{gnuest1}
g(\nu) \leq \nu\left[\kappa(a^*,1)-\kappa\left(a^*,1-\frac{a^*}{\nu}\right)\right].
\ee
Letting $\nu\to\infty$, we get from Lemma \ref{l:ka}(v) that
\be{gnest2}
\limsup_{\nu\to\infty} g(\nu) \leq a^* \left(\frac{\partial}{\partial b}
\kappa\right)(a^*,1) = \frac12\log\frac95.
\ee
Combine this with (\ref{fineq})(i) to get (\ref{fineq})(ii).
\epr

Proposition \ref{p:phtrinfchar} says that the free energy per step for an $AB$-block
exceeds that for an $AA$-block if and only the free energy per step for the single 
linear interface exceeds the free energy per step for an $AA$-block by a certain positive 
amount. This excess is needed to \emph{compensate} for the loss of entropy that occurs 
when the path runs along the interface for awhile before moving upwards from the 
interface to end at the diagonally opposite corner (recall Fig.\ 4). The constant 
$\frac12\log\frac95$ is special to our model.

\subsection{Criterion for $\psi_{AB}(\bar x)>\psi_{AA}(\bar x)$ and 
$\psi_{BA}(\bar y)>\psi_{BB}(\bar y)$}
\label{S2.4}

For all $(\alpha,\beta)\in\R^2$ and $a\geq 2$, we have
\be{psiineq}
\psi_{AB}(a) \geq \psi_{AA}(a), \qquad
\psi_{BA}(a) \geq \psi_{BB}(a).
\ee
The following gives a criterion for when strict inequality occurs and is the
analogue of Proposition \ref{p:phtrinfchar}.

\bp{p:psiineqchar}
For all $a\geq 2$, $\psi_{AB}(a)>\psi_{AA}(a)$ if and only if
\be{phinfcra}
\sup_{\mu \geq 1} \mu\left[\phi^\cI(\mu)-\frac12\alpha
-\frac12\log\left(\frac{a}{a-2}\right)\right]
>\frac12\log\left[\frac{4(a-2)(a-1)^2}{a}\right].
\ee
\ep

\bpr
Return to Lemma \ref{l:linkinf}. Fix $a\geq 2$. By (\ref{ka1psi}), (\ref{psiinflink})
and (\ref{feinfscal}), we have
\be{psidiffrel}
\begin{aligned}
&\psi_{AB}(a)-\psi_{AA}(a) = \sup_{0 \leq b \leq 1,\,a_1 \geq b,\,a_2 \geq 2-b,\,
a_1+a_2=a}\\
&\qquad\frac{a_1\phi^\cI(a_1/b)+a_2[\frac12\alpha+\kappa(a_2,1-b)]
-(a_1+a_2)\frac12\alpha-(a_1+a_2)\kappa(a_1+a_2,1)}
{a_1+a_2}.
\end{aligned}
\ee
The denominator is fixed. Put $\mu=a_1/b$ and rewrite the numerator as
\be{numrew}
\begin{aligned}
&\mu b\phi^\cI(\mu) + (a-\mu b)\left[\frac12\alpha+\kappa(a-\mu b,1-b)\right]
- a\left[\frac12\alpha+\kappa(a,1)\right]\\
&\qquad = \mu b\phi^\cI(\mu) - \mu b\frac12\alpha 
- [a\kappa(a,1)-(a-\mu b)\kappa(a-\mu b,1-b)].
\end{aligned}
\ee
By picking $s=a$, $t=1$, $u=a-\mu b$, $v=1-b$ in (\ref{chiconv}), we obtain that 
for all $\mu\geq 1$ and $0<b\leq 1$ with $(\mu-1)b\leq a-2$,
\be{kapprop}
a\kappa(a,1)-(a-\mu b)\kappa(a-\mu b,1-b)
> \mu b \left(\frac{\partial}{\partial a}(a\kappa)\right)(a,1) 
+ b \left(\frac{\partial}{\partial b}(a\kappa)\right)(a,1).
\ee
Since the right-hand side of (\ref{kapprop}) is $b$ times the derivative at $b=0$ 
of the left-hand side, it follows that the difference in (\ref{numrew}) is $\leq 0$ 
for all $0\leq b\leq 1$ if and only if its derivative at $b=0$ is $\leq 0$. This 
derivative equals
\be{diffb0}
\mu\phi^\cI(\mu)-\mu\frac12\alpha 
-\mu\left(\frac{\partial}{\partial a}(a\kappa)\right)(a,1)
-\left(\frac{\partial}{\partial b}(a\kappa)\right)(a,1).
\ee
After substituting the expressions for $\kappa(a,1)$, $(\frac{\partial}{\partial a}
\kappa)(a,1)$ and $(\frac{\partial}{\partial b}\kappa)(a,1)$ that we computed in Section 
\ref{S2.1.1} (recall (\ref{ka}), (\ref{maximisers}--\ref{maximisers*}), 
(\ref{parak}--\ref{parbk})), we find that (\ref{diffb0}) equals
\be{critb0}    
\mu\left[\phi^\cI(\mu)-\frac12\alpha
-\frac12\log\left(\frac{a}{a-2}\right)\right]
-\frac12\log\left[\frac{4(a-2)(a-1)^2}{a}\right].
\ee
Hence we get the claim.
\epr

For $\psi_{BA}(a)>\psi_{BB}(a)$ the same criterion applies as in (\ref{phinfcra}) with 
$\frac12\alpha$ replaced by $\frac12\beta$. (Recall that $\phi^\cI$ is symmetric in
$\alpha$ and $\beta$ by the remark made below (\ref{SBArel}).)

\subsection{Analysis of $F(\rho)$}
\label{S2.5}

In this section we analyse the variational problem in (\ref{Fdef}).

\bp{p:Fanal}
Let $(\alpha,\beta)\in\CONE$ and $\rho\in(0,1)$. Abbreviate $C=\alpha-\beta\geq 0$. 
The variational formula in {\rm (\ref{Fdef})} has unique maximisers $\bar x 
=\bar x(C,\rho)$ and $\bar y=\bar y(C,\rho)$ satisfying:\\
(i) $2 < \bar y < a^* < \bar x < \infty$ when $C>0$ and $\bar x=\bar y=a^*$ 
when $C=0$.\\
(ii) $u(\bar x)>v(\bar y)$ when $C>0$ and $u(\bar x) = v(\bar y)$ when $C=0$.\\
(iii) $\rho\mapsto\bar x(C,\rho)$ and $\rho\mapsto\bar y(C,\rho)$ are analytic 
and strictly decreasing on $(0,1)$ for all $C>0$.\\
(iv) $C\mapsto\bar x(C,\rho)$ and $C\mapsto\bar y(C,\rho)$ are analytic and 
strictly increasing, respectively, strictly decreasing on $(0,\infty)$
for all $\rho\in(0,1)$.\\
(v) As $\rho\ua 1$, $\bar x(C,\rho) \da a^*$ and $\bar y(C,\rho) \da 
10/(5-e^{-C})$ for all $C\geq 0$.\\
(vi) As $\rho\da 0$, $\bar x(C,\rho) \ua 10e^{-C}/(5e^{-C}-1)$ and 
$\bar y(C,\rho) \ua a^*$ when $0 \leq C <\log 5$, while $\bar x(C,\rho)
\ua\infty$ and $\bar y(C,\rho) \ua 2/(1-e^{-C})$ when $C \geq \log 5$.\\
(vii) As $C\ua\infty$, $\bar x(C,\rho)\ua\infty$ and $\bar y(C,\rho)\da 2$ 
for all $\rho\in (0,1)$.
\ep

\bpr
Fix $(\alpha,\beta)\in\CONE$ and $\rho\in (0,1)$. The supremum in (\ref{Fdef}) is 
attained at those $x,y$ that solve the equations
\be{xysol1}
\begin{aligned}
0 &= -\log 2 + \frac{(1-\rho)y}{2}(\alpha-\beta)+\frac{(1-\rho)y}{2}\log
\left(\frac{x(y-2)}{y(x-2)}\right)-\rho\log (x-2)-(1-\rho)\log (y-2),\\
0 &= -\log 2 + \frac{\rho x}{2}(\beta-\alpha)+\frac{\rho x}{2}\log
\left(\frac{y(x-2)}{x(y-2)}\right)-\rho\log(x-2)-(1-\rho)\log(y-2).
\end{aligned}
\ee
Multiplying the first relation by $\rho x$, the second relation by $(1-\rho)y$,
and adding them up, we get
\be{xysol2}
0 = -[\rho x+(1-\rho)y]\left\{\log 2+\rho\log(x-2)+(1-\rho)\log(y-2)\right\}.
\ee
Alternatively, subtracting the second relation from the first, we get
\be{xysol3}
0 = \frac12[\rho x+(1-\rho)y]\left\{(\alpha-\beta)
+\log\left(\frac{x(y-2)}{y(x-2)}\right)\right\}.
\ee
Hence, $x,y$ solve the equations
\be{xysol4}
\begin{aligned}
0 &= \log 2+\rho\log(x-2)+(1-\rho)\log(y-2),\\
0 &= (\alpha-\beta)+\log\left(\frac{x(y-2)}{y(x-2)}\right).
\end{aligned}
\ee  
These are two coupled equations depending on $\rho$, respectively, $C=\alpha-\beta$. 
Since the equations are linearly independent, their solution is unique.

\medskip\noindent
(i) Let $\bar x$ and $\bar y$ denote the unique solution of (\ref{xysol4}). Clearly, 
$\bar x = \bar y = a^*=\frac52$ when $C=0$. Suppose that $C>0$. Then it follows 
from the second line of (\ref{xysol4}) that $\bar x/(\bar x-2)<\bar y/(\bar y-2)$, 
or $\bar x>\bar y$. Moreover, it follows from the first line of (\ref{xysol4}) 
that it is not possible to have $\bar x>\bar y\geq a^*$ or $\bar y<\bar x\leq a^*$. 
Consequently,
\be{xysol5}
\bar y < a^* < \bar x.
\ee
The fact that $(\bar x,\bar y)\neq (\infty,2)$ follows from (\ref{xysol4}) as well.

\medskip\noindent
(ii) By (\ref{fgdef}),
\be{xysol7}
u(\bar x)-v(\bar y) = \frac{\alpha-\beta}{2}
+\left(\frac{1}{\bar x}-\frac{1}{\bar y}\right)\log 2 
+ \frac12\log\left(\frac{\bar x}{\bar y}\right) 
- \frac{\bar x-2}{2\bar x}\log(\bar x-2) 
+ \frac{\bar y-2}{2\bar y}\log(\bar y-2).
\ee 
Using (\ref{xysol4}), we may rewrite $\alpha-\beta$ and $\log 2$ in terms of
$\bar x,\bar y,\rho$. This gives, after some cancellations,
\be{xysol8}
u(\bar x)-v(\bar y) = \left(\frac{\rho}{\bar x}+\frac{1-\rho}{\bar y}\right)
\log\left(\frac{\bar x-2}{\bar y-2}\right).
\ee
This is $>0$ when $C>0$, because then $\bar x>\bar y$, and is $=0$ when $C=0$, 
because then $\bar x = \bar y$.

\medskip\noindent
(iii) The analyticity follows from the uniqueness of the solution of (\ref{xysol4})
and the implicit function theorem. From the second line of (\ref{xysol4}) it follows 
that $\rho\mapsto\bar x(C,\rho)$ and $\rho\mapsto\bar y(C,\rho)$ are either both 
non-increasing or both non-decreasing. Differentiating the first line of 
(\ref{xysol4}) w.r.t.\ $\rho$, we get 
\be{xysol9}
0 = \log\left(\frac{\bar x-2}{\bar y-2}\right)
+ \frac{\rho}{\bar x-2}\,\frac{\partial}{\partial\rho}\bar x 
+ \frac{1-\rho}{\bar y-2}\,\frac{\partial}{\partial\rho}\bar y.
\ee
Since $\bar x>\bar y$ when $C>0$, the sum of the last two terms is $<0$. Therefore 
it is not possible that $\frac{\partial}{\partial\rho}\bar x,\frac{\partial}{\partial\rho}
\bar y \geq 0$. Hence $\frac{\partial}{\partial\rho}\bar x,\frac{\partial}{\partial\rho}
\bar y < 0$.

\medskip\noindent
(iv) The analyticity again follows from the uniqueness of the solution of (\ref{xysol4})
and the implicit function theorem. From the first line of (\ref{xysol4}) it follows 
that $C\mapsto\bar x(C,\rho)$ and $C\mapsto \bar y(C,\rho)$ are either non-decreasing, 
respectively, non-increasing or vice versa. Differentiating the second line of 
(\ref{xysol4}) w.r.t.\ $C=\alpha-\beta$, we get 
\be{xysol10}
0 = 1 - \frac{2}{\bar x(\bar x-2)}\frac{\partial}{\partial C}\bar x
+ \frac{2}{\bar y(\bar y-2)}\frac{\partial}{\partial C}\bar y.
\ee
Since $\bar x,\bar y>2$, it is not possible that $\frac{\partial}{\partial C}\bar x
\leq 0 \leq \frac{\partial}{\partial C}\bar y$. Hence $\frac{\partial}{\partial C}
\bar y < 0 < \frac{\partial}{\partial C}\bar x$.

\medskip\noindent
(v) Abbreviate $\Delta=e^{-C} \in (0,1]$. Since $\bar y \leq a^*$, it follows from 
the first line of (\ref{xysol4}) that $\bar x \da a^*$ as $\rho\ua 1$. The second 
line of (\ref{xysol4}) therefore gives $\bar y/(\bar y-2) \ua 5\Delta^{-1}$, i.e.,
$\bar y \da 10/(5-\Delta)$.

\medskip\noindent
(vi) It follows from (\ref{xysol4}) that, as $\rho \da 0$, either $\bar x \ua A 
\in (a^*,\infty)$, $\bar y \ua a^*$ or $\bar x \ua \infty$, $\bar y \ua 2/(1-\Delta)$.
Since $\bar y \leq a^*=\frac52$, the latter is possible only when $\Delta\leq\frac15$. 
The former applies when $\Delta>\frac15$, in which case $A=10/(5-\Delta^{-1})
=10\Delta/(5\Delta-1)$.

\medskip\noindent
(vii)
This is immediate from (\ref{xysol4}).
\epr

\bl{Fanal}
$(\alpha,\beta)\mapsto F(\alpha,\beta;\rho)$ is analytic on $\R^2$
for all $\rho\in(0,1)$.
\el

\bpr
This is immediate from (\ref{fgdef}--\ref{Fdef}) and Proposition \ref{p:Fanal}(iv). 
\epr


\section{Free energy of the polymer}
\label{S3}

In Section \ref{S3.1} we prove Theorem \ref{feiden}. In Section \ref{S3.2}
we analyse the set $\cR(p)$ in (\ref{Rdef}) and the supremum $\rho^*(p)$ in 
(\ref{rho*def}).

\subsection{Proof of Theorem \ref{feiden}}
\label{S3.1}

(i) Write out the partition sum in (\ref{fedef}) in terms of partition
sums in successive blocks:
\be{partrew}
\begin{aligned}
Z^{\omega,\Omega}_{n,L_n} &= \sum_{N=1}^{n/2L_n}\, \sum_{(\Pi_i)_{i=1}^N}\,\,
\sum_{u_1=2L_n}^\infty \cdots \sum_{u_N=2L_n}^\infty
\Big[\prod_{i=1}^{N-1} 
Z^{\sigma^{u_1+ \cdots + u_{i-1}}\omega}_{u_i}(t^{\Omega}(\Pi_i))\Big]\\
&\qquad \times Z^{\sigma^{u_1 + \cdots + u_{N-1}}\omega}_{n-(u_1+ \cdots + u_{N-1})}
((t^{\Omega}(\Pi_i))\,\,
1\{u_1+ \cdots + u_{N-1} \leq n < u_1 + \cdots + u_{N-1} + u_N\}.
\end{aligned}
\ee
Here, $N$ counts the number of blocks traversed, $\sigma$ is the left-shift
acting on $\omega$, $\Pi_i$ is the $i$-th step of the coarse-grained path, 
$u_i$ counts the number of steps spent in the $i$-th block diagonally traversed 
by $\Pi_i$, $t^{\Omega}(\Pi_i)$ labels the type of the $i$-th block in $\Omega$, 
and $Z^{\omega}_u(t)$ is the partition sum for spending $u$ steps in a block of 
type $t$. We want to derive the asymptotics of this expression as $n\to\infty$. 
\emph{For reasons of space the argument below is somewhat sketchy, but the 
technical details are easy to fill in.}

First, for the computation we pretend that after $n$ steps the path has just 
completed traversing a block, i.e., we replace the indicator in (\ref{partrew})
by
\be{indic}
1\{u_1 + \cdots + u_N=n\}.
\ee
The error made in doing so is at most a factor $O(n)$. Next, in (\ref{partrew}) 
we insert a weight $2^{-N}$ under the sum over $(\Pi_i)_{i=1}^N$, scale $u_i$ 
by putting $v_i=u_i/L_n$, change the sum over $u_i$ to an integral over $v_i$, 
and insert a weight $e^{-(v_i-2)}$ under the integral. This gives
\be{intrew}
\begin{aligned}
Z^{\omega,\Omega}_{n,L_n} &= O(n) \sum_{N=1}^{n/2L_n} 
\left(2^N e^{(n/L_n)-2N} L_n^N\right)\\
&\qquad\times 
\Big[\sum_{(\Pi_i)_{i=1}^N} 2^{-N}\,\,
\int_2^\infty dv_1\,e^{-(v_1-2)} \cdots 
\int_2^\infty dv_N\,e^{-(v_N-2)}\\
&\qquad\qquad\times \prod_{i=1}^N 
Z^{\sigma^{(v_1+\cdots v_{i-1})L_n}\omega}_{v_iL_n}(t^{\Omega}(\Pi_i))\,\,
1\{v_1 + \cdots + v_N = n/L_n\}\Big],
\end{aligned}
\ee
where the term between round brackets compensates for the insertion of
the weights (and can be computed because of (\ref{indic})), while roundoff 
errors (coming from turning sums into integrals) disappear into the error term. 
The factor between round brackets is $e^{o(n)}$ and therefore is negligible. 

The point of the rewrite in (\ref{intrew}) is that the sum over $(\Pi_i)_{i=1}^N$
and the integrals over $v_i$ are \emph{normalised}. Therefore we can now introduce
two independent sequences of \emph{random variables}, 
\be{rvdef}
\hat\Pi = \{\hat \Pi_i\colon\,i \in \N\}, \qquad 
\hat v = \{\hat v_i\colon\,i \in \N\},
\ee
which describe a random uniform coarse-grained path, repectively,
a random sequence of scaled times that are i.i.d.\ and ${\rm Exp}(1)$ 
distributed on $[2,\infty)$. In terms of these random variables we can 
rewrite (\ref{intrew}) as
\be{intrew*}
Z^{\omega,\Omega}_{n,L_n} = e^{o(n)} \sum_{N=1}^{n/2L_n}\\
\Big\langle \prod_{i=1}^N 
Z^{\sigma^{(\hat v_1+\cdots \hat v_{i-1})L_n}\omega}_{\hat v_iL_n}
(t^{\Omega}(\hat \Pi_i))\,\,
1\{\hat v_1 + \cdots + \hat v_N = n/L_n\}\Big]\Big\rangle,
\ee
where $\langle\cdot\rangle$ denotes expectation with respect to $(\hat\Pi,\hat v)$. 
As $L_n\to\infty$ we have, for every fixed realisation of $(\hat\Pi,\hat v)$, 
\be{frlimrel}
Z^{\sigma^{(\hat v_1+\cdots+\hat v_{i-1})L_n}\omega}_{\hat v_iL_n}
(t^{\Omega}(\hat\Pi_i)) = \exp\left\{L_n [1+o(1)]\,\hat v_i\,
\psi_{t^{\Omega}(\hat\Pi_i)}(\hat v_i)\right\} 
\qquad \omega-a.s.
\ee
with $\psi_{kl}(\hat v_i)$ the free energy per step in a $kl$-block where 
the path spends $\hat v_iL_n$ steps, defined in (\ref{fekl}). Here we use 
that the distribution of $\omega$ is invariant under shifts and that
$(\hat v_1+\cdots+\hat v_{i-1})L_n$ is independent of $\omega$. 

Because of (\ref{intrew*}) and (\ref{frlimrel}), we are in a position to use 
large deviation theory (for background see e.g.\ den Hollander \cite{dH00},
Chapters I and II). To that end, we introduce the \emph{empirical distribution}
\be{empdef}
\cE_N^\Omega = \cE_N^\Omega(\hat\Pi,\hat v) 
= \frac1N \sum_{i=1}^N \delta_{(t^{\Omega}(\hat \Pi_i),\hat v_i)},
\ee
where $\delta_{(t,v)}$ is the unit measure at $(t,v)$. This $\cE_N^\Omega$
counts the frequency at which the $(t^{\Omega}(\hat \Pi_i),\hat v_i)$
assume values in the space
\be{Thetadef}
\Theta = \{AA,AB,BA,BB\} \times [2,\infty)
\ee
and is an element of $\cP(\Theta)$, the set of probability distributions $\Theta$.
With the help of (\ref{frlimrel}), we may rewrite (\ref{intrew*}) as
\be{intrewemp}
Z^{\omega,\Omega}_{n,L_n} = e^{o(n)} \sum_{N=1}^{n/2L_n}
\Big\langle \exp\Big\{NL_n[1+o(1)]\,(\cE_N^\Omega,h_1)_*\Big\}
\,\,1\{(\cE_N^\Omega,h_2)_*=n/NL_n\}
\Big\rangle,
\ee
where we introduce two functions on $\Theta$,
\be{hdefs}
h_1(t,v) = v\psi_t(v) \quad \hbox{and} \quad h_2(t,v) = v \quad \hbox{for} \quad
(t,v) \in \Theta,
\ee
and put $(\mu,h)_*=\int_\Theta hd\mu$ for $\mu\in\cP(\Theta)$. At this point,
our partition sum has been rewritten in terms of an expectation w.r.t.\
the empirical distribution $\cE_N^\Omega$ (which depends on $(\hat\Pi,\hat v)$). 

Our next step is to scale $N$ by putting $M=NL_n/n$ and to rewrite (\ref{intrewemp}) 
as
\be{intrewfin}
Z^{\omega,\Omega}_{n,L_n} = e^{o(n)} \int_0^{\frac12} dM\,
\left\langle\exp\left[n 
\frac{(\cE_N^\Omega,h_1)_*}{(\cE_N^\Omega,h_2)_*}
\right]\,1\left\{(\cE_N^\Omega,h_2)_*=\frac1M\right\}\right\rangle,
\ee
where again roundoff errors disappear into the error term. The idea that we now use 
is that $(\cE^{\Omega}_N)_{N\in\N}$ satisfies the \emph{large deviation principle} on 
$\cP(\Theta)$ with rate $N$ and with some rate function $\mu \mapsto I(\mu)$ that 
has compact level sets. This rate function will be $\Omega$-a.s.\ constant. For 
techniques on how to prove this, we refer to Comets \cite{C89}, Greven and den 
Hollander \cite{GH92} and Sepp\"al\"ainen \cite{Se95}. Thus, the probability of 
$\cE^{\Omega}_N$ being close to some $\mu\in\cP(\Theta)$ is $\exp\{-N[1+o(1)]I(\mu)\}$. 
Since $N=Mn/L_n=o(n)$, we conclude that these large deviations have a \emph{negligible 
cost}. Using (\ref{intrewfin}), together with the relation
\be{intout}
\int_0^{\frac12} dM\,1\left\{(\cE_N^\Omega,h_2)_*=\frac1M\right\} = 
1/O(n)
\ee
(which follows from reversing the calculations above), we arrive at
\be{felimres}  
f = \lim_{n\to\infty} \frac1n \log Z^{\omega,\Omega}_{n,L_n}
= \sup_{\mu\in\cP(\Theta)} \frac{(\mu,h_1)_*}{(\mu,h_2)_*} \qquad \omega-a.s.
\ee
But the supremum in the right-hand side is precisely the formula for $f$ 
stated in Theorem \ref{feiden}(i). The free energy is trivially finite 
(recall (\ref{Hamiltonian}), (\ref{fedef}) and (\ref{sa})).

\medskip\noindent
{\bf Remark:} In the above, somewhat sketchy, computation the introduction of 
the rate function can be avoided by appealing to concentration of measure 
estimates, which are a crude yet flexible form of large deviations (see Madras and 
Whittington \cite{MW02} for an application of this technique in the context 
of self-avoiding random copolymers near a single linear interface).

\medskip\noindent
(ii) The proof is elementary. Fix $\omega$ and $\Omega$, and rewrite the partition sum
in (\ref{fedef}) as
\be{ferew}
Z^{\omega,\Omega}_{n,L_n} = \sum_{v_A,v_B} c^{\omega,\Omega}_{n,L_n}(v_A,v_B)
\,\,e^{\alpha v_A + \beta v_B}
\ee
with
\be{vABdef}
\begin{aligned}
&c^{\omega,\Omega}_{n,L_n}(v_A,v_B)\\ 
&\qquad = \left|\left\{
\pi\in\cW_{n,L_n}\colon\,
\sum_{i=1}^n 1\{\omega_i=\Omega_{\pi_i}=A\}=v_A,\,
\sum_{i=1}^n 1\{\omega_i=\Omega_{\pi_i}=B\}=v_B
\right\}\right|.
\end{aligned}
\ee
Pick any $\alpha_1,\alpha_2,\beta_1,\beta_2\in\R$. Then, by Cauchy-Schwarz applied to
(\ref{ferew}),
\be{CS}
Z^{\omega,\Omega}_{n,L_n}\left(\frac{\alpha_1+\alpha_2}{2},\frac{\beta_1+\beta_2}{2}\right) 
\leq Z^{\omega,\Omega}_{n,L_n}(\alpha_1,\beta_1)Z^{\omega,\Omega}_{n,L_n}(\alpha_2,\beta_2), 
\ee
from which the claim follows (recall (\ref{fedef}) and (\ref{sa})).

\medskip\noindent
(iii) According to Proposition \ref{Rpprop}(i) below, $p\mapsto\cR(p)$ is continuous
in the Hausdorff metric. It is therefore immediate from the variational representation 
of $f$ in (\ref{fevar}), together with the continuity of $\psi_{kl}$ stated in Lemma
\ref{psiprop}(i), that $p\mapsto f(\alpha,\beta;p)$ is continuous.  

\subsection{Analysis of $\cR(p)$ and $\rho^*(p)$}
\label{S3.2}

The following proposition is crucial for the analysis of the phase transition
curves. Recall (\ref{Rdef}), (\ref{rho*def}) and Fig.\ 6. The elements of 
$\cR(p)$ are matrices
\be{matr}
\left(\ba{ll} \rho_{AA} &\rho_{AB}\\ \rho_{BA} &\rho_{BB} \ea\right)
\ee
whose elements are non-negative and sum up to 1.

\bp{Rpprop}
(i) $p \mapsto \cR(p)$ is continuous in the Hausdorff metric.\\
(ii) If $p \geq p_c$, then
\be{Rpa}
\left(\ba{ll} 1-\gamma &\gamma\\ 0 &0 \ea\right)
\in \cR(p) \quad \hbox{ for some } \gamma \in (0,1).
\ee
(iii) If $p<p_c$, then
\be{Rpb}
\left(\ba{ll} 1-\gamma &\gamma\\ 0 &0 \ea\right)
\notin \cR(p) \quad \hbox{ for all } \gamma \in [0,1].
\ee
\ep

\bpr
(i) Return to (\ref{fracdef}--\ref{Rdef}). Pick $0<p<p'<1$. Let $\Omega$ and
$\Omega'$ be two typical percolation configurations with parameter $p$ and 
$p'$, respectively, coupled such that the set of $A$'s in $\Omega'$ contains
the set of $A$'s in $\Omega$. Then 
\be{comprho}
\limsup_{n\to\infty} |\rho^\Omega_{kl}(\Pi,n)-\rho^{\Omega'}_{kl}(\Pi,n)| 
\leq 2(p'-p) \quad \Omega,\Omega'-a.s. 
\qquad \forall\,k,l\in\{A,B\},\,\Pi\in\cW.  
\ee
Therefore, for all $k,l\in\{A,B\}$, we have $|(\cR^\Omega-\cR^{\Omega'})_{kl}|
\leq 2(p'-p)$ $\Omega,\Omega'$-a.s.\ and hence $|[\cR(p)-\cR(p')]_{kl}|\leq 
2(p'-p)$.

\noindent
(ii) If $p>p_c$, then there is an infinite cluster of $A$-blocks with
a strictly positive density. For the coarse-grained path $\Pi$ that
moves to the infinite cluster and afterwards stays inside this
cluster, we have $\rho_A=\rho_{AA}+\rho_{AB}=1$. For the coarse-grained 
path $\Pi$ that moves to the infinite cluster, afterwards stays inside
this cluster, but follows its boundary as much as possible, we have 
$0<\rho_{AB}=1-\rho_{AA}<1$, which proves the claim. Since $p\mapsto
\cR(p)$ is continuous in the Hausdorff metric, and since lowering $p$ 
increases the density of the $B$ blocks, the claim trivially extends to 
$p=p_c$.

\noindent
(iii) If $p<p_c$, then there is no infinite cluster of $A$-blocks. In fact,
any coarse-grained path $\Pi$ visits $B$-blocks with a strictly positive 
density (as follows from a coupling argument similar as in (i)). Hence,
$\rho_{A}=\rho_{AA}+\rho_{AB}<1$ for all $\Pi$, and in fact $\sup_{\Pi\in\cW}
\rho_A(\Pi)<1$. Since $\cR(p)$ is a closed set, this proves the claim.
\epr

It follows from Proposition \ref{Rpprop} and the fact that $\cR(p)$ is closed
that $\rho^*(p)$ defined in (\ref{rho*def}) has the qualitative properties 
indicated in Fig.\ 6: $p\mapsto\rho^*(p)$ is continuous and non-decreasing 
on $(0,1)$, $\rho^*(p)=1$ for $p\in [p_c,1)$ and $\rho^*(p)\in (0,1)$ for 
$p\in (0,p_c)$.


\section{Analysis of the critical curve}
\label{S4}

In Section \ref{S4.1} we prove Theorems \ref{fesupcr}, \ref{phtriden} and 
\ref{phtrcurve} for $p \geq p_c$. In Section \ref{S4.2} we prove Theorems 
\ref{fesubcr}, \ref{phtriden*} and \ref{phtrcurve*} for $p < p_c$. In
Section \ref{S4.3} we make some observations about the separation of
the localised phase into two subphases for $p<p_c$.

\subsection{Supercritical case $p \geq p_c$}
\label{S4.1}

\subsubsection{Proof of Theorems \ref{fesupcr} and \ref{phtriden}}
\label{S4.1.1}

The following proposition proves Theorems \ref{fesupcr} and \ref{phtriden}.
Recall from (\ref{Sineq}) that $S_{AB} \geq S_{AA}$.

\bp{phtrcond}
Fix $p \geq p_c$.\\
(i) If $S_{AB}=S_{AA}$, then $f=S_{AA}$.\\
(ii) If $S_{AB}>S_{AA}$, then $f>S_{AA}$.
\ep

\bpr
The proof uses Theorem \ref{feiden}(i) and Proposition \ref{Rpprop}(ii), in 
combination with the inequalities
\be{Sineqsuppl}
S_{BB} \leq S_{AA}, \qquad S_{BA} \leq S_{AB},
\ee
which hold because $\beta\leq\alpha$ (recall Propositions \ref{p:fe} and \ref{p:SAB}).

\noindent
(i) Suppose that $S_{AB}=S_{AA}$. Then, because $\psi_{kl}(a) \leq S_{kl}$ for all
$a\geq 2$ and $k,l\in\{A,B\}$ by (\ref{Sdef}), Theorem \ref{feiden}(i) and 
(\ref{Sineqsuppl}) yield
\be{festa}
f \leq \sup_{(a_{kl})\in\cA}\,\, \sup_{(\rho_{kl}) \in \cR(p)}
\frac{\sum_{k,l} \rho_{kl} a_{kl} S_{kl}}
{\sum_{k,l} \rho_{kl} a_{kl}}
\leq \sup_{k,l} S_{kl} = S_{AB} = S_{AA}.
\ee
On the other hand, Theorem \ref{feiden}(i) and Proposition \ref{Rpprop}(ii) yield
\be{festc*}
f \geq \frac{(1-\gamma)\bar a_{AA}S_{AA} + \gamma \bar a_{AB}S_{AB}}
{(1-\gamma) \bar a_{AA} + \gamma \bar a_{AB}}
= S_{AA},
\ee
where $\bar a_{AA},\bar a_{AB}$ are the maximisers of $S_{AA},S_{AB}$ (and the 
value of $\gamma$ is irrelevant). Combine (\ref{festa}) and (\ref{festc*}) to get
$f=S_{AA}$.

\noindent
(ii) Suppose that $S_{AB}>S_{AA}$. Then Theorem \ref{feiden}(i) and Proposition 
\ref{Rpprop}(ii) with $0<\gamma<1$ yield
\be{festc}
f \geq \frac{(1-\gamma)\bar a_{AA}S_{AA} + \gamma \bar a_{AB}S_{AB}}
{(1-\gamma) \bar a_{AA} + \gamma \bar a_{AB}}
> S_{AA}
\ee
(here it is important that $\gamma>0$). 
\epr

We see from Proposition \ref{phtrcond} that $\cD$ corresponds to the situation
where the polymer is fully $A$-delocalized ($f=S_{AA}$), while $\cL$ corresponds
to the situation where the polymer is partially $AB$-delocalized ($S_{AA}<f<S_{AB}$).  

\subsubsection{Proof of Theorem \ref{phtrcurve}}
\label{S4.1.2}

Since $S_{AA}(\alpha,\beta)$ does not depend on $\beta$, and $\beta\mapsto 
S_{AB}(\alpha,\beta)$ is continuous and non-decreasing on $\R$ for every 
$\alpha\in\R$ by Lemma \ref{psiprop}(ii), the boundary between $\cD$ and 
$\cL$ is a continuous function in $\CONE$. We denote this function by 
$\alpha\mapsto\beta_c(\alpha)$. 

We first show that the curve is concave. To that end, pick any $\alpha_1<\alpha_2$, 
and consider the points $(\alpha_1,\beta_c(\alpha_1))$ and $(\alpha_2,\beta_c(\alpha_2))$
on the curve. Let $(\alpha_3,\beta_3)$ be the midpoint of the line connecting the two. 
We want to show that $\beta_3 \leq \beta_c(\alpha_3)$. By the convexity of $f$, stated 
in Theorem \ref{feiden}(ii), we have
\be{midineq}
f(\alpha_3,\beta_3) \leq
\frac12[f(\alpha_1,\beta_c(\alpha_1))+f(\alpha_2,\beta_c(\alpha_2))].
\ee
Since the curve itself is part of $\cD$ (recall (\ref{DLdef})), it follows from 
Propositions \ref{p:fe} and \ref{phtrcond}(i) that the right-hand side of (\ref{midineq}) 
equals
\be{eqf}
\frac12\left[\left(\frac12\alpha_1+\kappa\right)
+\left(\frac12\alpha_2+\kappa\right)\right]
= \frac12\left(\frac{\alpha_1+\alpha_2}{2}\right)+\kappa
= \frac12\alpha_3+\kappa.
\ee  
Thus, $f(\alpha_3,\beta_3) \leq \frac12\alpha_3 + \kappa$. But, by Proposition 
\ref{phtrcond}, the reverse inequality is true always, and so equality holds. 
Consequently, $(\alpha_3,\beta_3)\in\cD$, which proves the claim that $\beta_3
\leq\beta_c(\alpha_3)$.

The concavity in combination with the lower bound in part (i) of the following 
lemma show that the curve is non-decreasing.

The following lemma settles most of Theorem \ref{phtrcurve}.

\bl{phtrest}
Fix $p \geq p_c$.\\ 
(i) $\beta_c(\alpha) \geq \log (2-e^{-\alpha})$ for all $\alpha\geq 0$.\\
(ii) $\beta_c(\alpha) < 8\log 3$ for all $\alpha\geq 0$.\\
(iii) $\beta_c(\alpha)=\alpha$ for all $0\leq\alpha\leq\alpha_0$, where 
$\alpha_0$ is the number defined in the first line of {\rm (\ref{alpha01})}.
\el

\bpr
(i) We have, recalling (\ref{feinf}--\ref{Zinf}),
\be{phiIrel}
\phi^\cI(\mu) = \lim_{L\to\infty} \frac{1}{L} \log Z^{\omega,\cI}_L(\mu)
\quad \omega-a.s.
\ee
with
\be{ZIrel}
\begin{aligned}
Z^{\omega,\cI}_L(\mu) &= \sum_{\pi\in\cW_{\mu L,L}} 
\exp\left[-H^{\omega,\cI}_{\mu L}(\pi)\right]\\
H^{\omega,\cI}_{\mu L}(\pi) 
&= - \sum_{i=1}^{\mu L}\Big(\alpha 1\{\omega_i=A,\pi_i>0\}
+\beta 1\{\omega_i=B,\pi_i \leq 0\}\Big).
\end{aligned}
\ee
We will derive an upper bound on $\phi^\cI(\mu)$ by doing a so-called first-order 
annealed estimate (also referred to as a first-order Morita approximation; see
Orlandini \emph{et al} \cite{ORW02}). This estimate consists in writing
\be{Morita1}
H^{\omega,\cI}_{\mu L}(\pi) 
= -\sum_{i=1}^{\mu L} \alpha 1\{\omega_i=A\}
- \sum_{i=1}^{\mu L} 1\{\pi_i \leq 0\}\,
\Big(- \alpha 1\{\omega_i=A\}+\beta 1\{\omega_i=B\}\Big),
\ee
using that the first term is $-\mu L\frac12\alpha [1+o(1)]$ $\omega$-a.s.\
as $L\to\infty$ and is independent of $\pi$, substituting this into (\ref{ZIrel}) 
and performing an expectation over $\omega$. This gives 
\be{Morita2}
\begin{aligned}
&\left\langle\log Z^{\omega,\cI}_L(\mu)\right\rangle\\
&\qquad \leq \mu L \frac12\alpha[1+o(1)]
+ \log \sum_{\pi\in\cW_{\mu L,L}}
\prod_{i=1}^{\mu L} 1\{\pi_i\leq 0\}\,
\left\langle e^{-\alpha 1\{\omega_i=A\}+\beta 1\{\omega_i=B\}}\right\rangle\\
&\qquad \leq \mu L \frac12\alpha[1+o(1)] + \mu L[\hat\kappa(\mu)+o(1)]
+\mu L \log\left(\frac12 e^{-\alpha}+\frac12 e^\beta\right),
\end{aligned}
\ee
where $\langle\cdot\rangle$ denotes expectation over $\omega$, we use Jensen's 
inequality as well as the i.i.d.\ property of $\omega$ in the first line,
and we use Lemma \ref{l:kamu}(i) in the second line. Consequently,
\be{Morita3}
\phi^\cI(\mu) =\lim_{L\to\infty} \frac{1}{\mu L} 
\left\langle\log Z^{\omega,\cI}_L(\mu)\right\rangle\\
\leq \frac12\alpha + \hat\kappa(\mu) + 
\log\left(\frac12 e^{-\alpha}+\frac12 e^\beta\right). 
\ee

Suppose that
\be{Morita4}
\log\left(\frac12 e^{-\alpha}+\frac12 e^\beta\right) \leq 0.
\ee
Then substitution of (\ref{Morita3}) into (\ref{diffrew}) gives
\be{Morita5}
S_{AB}-S_{AA} \leq \sup_{\mu\geq 1,\,\nu\geq 1}
\frac{\mu[\hat\kappa(\mu)-\frac12 \log 5]-\nu[\frac12 \log 5 - f(\nu)]}{\mu+\nu}.
\ee
But the right-hand side is the same as $S_{AB}-S_{AA}$ when $\alpha=\beta=0$
(as can be seen from (\ref{diffrew}) because $\phi^\cI(\mu)=\hat\kappa(\mu)$ when
$\alpha=\beta=0$), and therefore is equal to 0. Hence, recalling (\ref{Ldef}), we 
find that (\ref{Morita4}) implies that $(\alpha,\beta)\in\cD$. Consequently,
\be{Morita6}
\log\left(\frac12 e^{-\alpha}+\frac12 e^{\beta_c(\alpha)}\right) \geq 0
\qquad \hbox{for all } \alpha\geq 0,
\ee
which gives the lower bound that is claimed.

\medskip\noindent
(ii) We will show that there exists a $\mu_0>1$ such that
\be{mucond1}
\mu_0[\phi^{\cI}(\mu_0)-S_{AA}] > \frac12 \log \frac95
\quad \mbox{ for all } \alpha\geq 0 \hbox{ when } \beta \geq 8\log 3.
\ee
This will prove the claim via Theorem \ref{phtriden} and Proposition 
\ref{p:phtrinfchar}.

Consider the polymer along the single infinite interface $\cI$ introduced
in Section \ref{S2.2.2}. Fix $\omega$. In $\omega$, look for the strings
of $B$'s that are followed by a string of at least three $A$'s. Call these 
$B$-strings ``good'', and call all other $B$-strings ``bad''. Let $\pi(\omega)$ 
be the path that starts at $(0,0)$, steps to $(0,1)$ and proceeds as follows. 
Each time a good $B$-string comes up, the path moves down from height 1 to 
height 0 during the step that carries the $A$ just preceeding the good $B$-string, 
moves at height 0 during the steps that carry the $B$'s inside the string, 
moves up from height 0 to height 1 during the step that carries the first $A$ 
after the string, and moves at height 1 during the step that carries the second 
$A$ after the string. The third $A$ can be used to either move from height 1 
to height 0 in case the next good $B$-string comes up immediately, or to 
move at height 0 in case it is not. When a bad $B$-string comes up, the 
path stays at height 1.   

Along $\pi(\omega)$, we have that all the $A$'s lie in the upper halfplane,
all the bad $B$-strings lie in the upper halfplane, while all the good 
$B$-strings lie in the interface. Asymptotically, the good $B$-strings 
contain $\frac14$-th of the $B$'s. Hence $\frac18$-th of the steps carry 
a $B$ that is in a good $B$-string. Moreover, the number of steps between 
heights 0 and 1 is $\frac18$ times the number of steps at heights 0 and 1 
(the average length of a good $B$-string is 2), and so $\pi(\omega)$ travels 
a distance $L$ in time $\frac98 L$ for $L$ large, which corresponds to 
$\mu=\mu_0=\frac98$. Thus, the contribution of $\pi(\omega)$ to the Hamiltonian 
in (\ref{Zinf}) equals
\be{Hcontr}
H^{\omega,\cI}_{\mu_0L}(\pi(\omega)) 
= - \mu_0L \left(\frac12\alpha + \frac18\beta\right)[1+o(1)] \qquad \omega-a.s.
\mbox{ as } L\to\infty.
\ee
Therefore, recalling (\ref{fesainf}) and (\ref{feinfscal}), we have
\be{psicontr}
\phi^{\cI}(\mu_0) \geq \frac12\alpha + \frac18\beta.
\ee
Via (\ref{fe}) this gives
\be{diffest}
\mu_0[\phi^{\cI}(\mu_0)-S_{AA}] 
\geq \mu_0 \left(\frac18\beta-\frac12 \log 5\right).
\ee
Consequently, the inequality in (\ref{mucond1}) holds as soon as
\be{betaineq}
\frac18\beta > \frac89 \log 3 + \frac{1}{18} \log 5.
\ee
Since the right-hand side is strictly smaller than $\log 3$, this proves
the claim. 

\medskip\noindent
(iii) Pick $\alpha=\beta$. Then, by (\ref{phirelbds}), $\phi^I(\mu)\leq\alpha
+\hat\kappa(\mu)$. If $\alpha \in [0,\alpha_0]$, with $\alpha_0$ given by the 
first line of (\ref{alpha01}), then this bound in combination with Propositions 
\ref{p:fe} and \ref{p:phtrinfchar} gives $S_{AB}=S_{AA}$. Thus, $\{(\alpha,\alpha)
\colon\,\alpha\in [0,\alpha_0]\} \subset \cD$. 
\epr

Lemma \ref{phtrest}, together with the concavity of $\alpha\mapsto\beta_c(\alpha)$
(shown prior to Lemma \ref{phtrest}), proves Theorem \ref{phtrcurve}, except for 
the slope discontinuity stated in (\ref{slopedis}). But the latter follows from 
the fact that if the piece of $\beta_c$ on $[\alpha^*,\infty)$ is analytically 
continued outside $\CONE$, then it hits the vertical axis at a strictly positive 
value, namely, $\alpha_0$ defined in the first line of (\ref{alpha01}). Indeed, 
for $\alpha=0$ we have $\phi^{\cI}(\mu) = \frac12\beta+\hat\kappa(\mu)$, because 
there is zero exponential cost for the path to stay in the lower halfplane 
(recall (\ref{feinf}--\ref{fesainf}) and (\ref{feinfscal})). Consequently, the 
criterion for delocalization in Proposition \ref{p:phtrinfchar}, $S_{AB}=S_{AA}$, 
reduces to
\be{redcrit}
\sup_{\mu\geq 1} \mu\left[\hat\kappa(\mu)+\frac12\beta-\frac12\log 5\right]
\leq \frac12\log\frac95.
\ee    
This is true precisely when $\beta\leq\alpha_0$.

\subsection{Subcritical case $p < p_c$}
\label{S4.2}

\subsubsection{Proof of Theorems \ref{fesubcr} and \ref{phtriden*}}
\label{S4.2.1}

\bp{p:fred}
If $\psi_{AB}\equiv\psi_{AA}$ and $\psi_{BA}\equiv\psi_{BB}$, then
$f=F(\rho^*(p))$.
\ep

\bpr
Suppose that $(\alpha,\beta)$ is such that $\psi_{AB}(a)=\psi_{AA}(a)$ and 
$\psi_{BA}(a)=\psi_{BB}(a)$ for all $a\geq 2$. Then the variational formula in
(\ref{fevar}) reduces to
\be{fvarred1}
\begin{aligned}
&f = \sup_{(a_{kl})\in\cA}\,\sup_{(\rho_{kl})\in\cR(p)}\\
&\frac{
\rho_{AA}a_{AA}\psi_{AA}(a_{AA})+\rho_{AB}a_{AB}\psi_{AA}(a_{AB})
+\rho_{BA}a_{BA}\psi_{BB}(a_{BA})+\rho_{BB}a_{BB}\psi_{BB}(a_{BB})
}{
\rho_{AA}a_{AA}+\rho_{AB}a_{AB}+\rho_{BA}a_{BA}+\rho_{BB}a_{BB}
}.
\end{aligned}
\ee
Define
\be{fvarred3}
\tilde a_{AA} = \frac{\rho_{AA}}{\rho_A}a_{AA}
+ \frac{\rho_{AB}}{\rho_A}a_{AB},
\qquad
\tilde a_{BB} = \frac{\rho_{BA}}{\rho_B}a_{BA}
+ \frac{\rho_{BB}}{\rho_B}a_{BB},
\ee
where $\rho_A=\rho_{AA}+\rho_{AB}$ and $\rho_B=\rho_{BA}+\rho_{BB}$. By Lemma 
\ref{l:ka}(ii) and (\ref{ka1psi}), $a\mapsto a\psi_{AA}(a)$ and $a\mapsto 
a\psi_{BB}(a)$ are concave on $[2,\infty)$. Hence, the numerator in (\ref{fvarred1}) 
can be bounded above as
\be{fvarred2}
\begin{aligned}
&\rho_A\left[\frac{\rho_{AA}}{\rho_A}a_{AA}\psi_{AA}(a_{AA})
+ \frac{\rho_{AB}}{\rho_A}a_{AB}\psi_{AA}(a_{AB})\right]\\
&\qquad\qquad + \rho_B\left[\frac{\rho_{BA}}{\rho_B}a_{BA}\psi_{BB}(a_{BA})
+ \frac{\rho_{BB}}{\rho_B}a_{BB}\psi_{BB}(a_{BB})\right]\\
&\qquad \leq \rho_A\tilde a_{AA}\psi_{AA}(\tilde a_{AA}) 
+ \rho_B\tilde a_{BB}\psi_{BB}(\tilde a_{BB}),
\end{aligned} 
\ee
while the denominator in (\ref{fvarred1}) equals
\be{fvarden}
\rho_A\tilde a_{AA} + \rho_B\tilde a_{BB}.
\ee
For any choice of $(\rho_{kl})\in\cR(p)$, as $(a_{kl})$ runs through $\cA$, 
$\tilde a_{AA}$ and $\tilde a_{BB}$ run through all values $\geq 2$. Moreover,
equality can be achieved in (\ref{fvarred2}) by picking $a_{AA}=a_{AB}$ and
$a_{BA}=a_{BB}$. Hence (\ref{fvarred1}) reduces to
\be{fvarred4}
f = \sup_{(\rho_{kl})\in\cR(p)} F(\rho_{A})
\ee
with $F(\rho)$ given by (\ref{Fdef}). Thus, it remains to show that the
supremum is taken at $\rho^*(p)=\sup_{(\rho_{kl})\in\cR(p)} \rho_A$. 

For $\rho\in(0,1)$, let $\bar x=\bar x(\rho)$ and $\bar y=\bar y(\rho)$ denote
the unique maximers of (\ref{Fdef}). Then, for any $\rho_1,\rho_2\in(0,1)$,
\be{fvarred5}
F(\rho_1)-F(\rho_2) = R(\rho_1,\bar x(\rho_1),\bar y(\rho_1))
- R(\rho_2,\bar x(\rho_2),\bar y(\rho_2))
\ee
with
\be{fvarrred6}
R(\rho,x,y) = \frac{\rho xu(x)-(1-\rho) yu(y)}{\rho x+(1-\rho)y}.
\ee
Since $R(\rho_1,\bar x(\rho_1),\bar y(\rho_1))\geq R(\rho_1,\bar x(\rho_2),\bar y(\rho_2))$,
we have
\be{fvarrred7}
\frac{\partial}{\partial\rho}F(\rho)
\geq \left(\frac{\partial}{\partial\rho}R\right)(\rho,\bar x(\rho),\bar y(\rho)).
\ee
Now compute
\be{fvarred8}
\left(\frac{\partial}{\partial\rho}R\right)(\rho,x,y)
= \frac{xy[u(x)-v(y)]}{[\rho x+(1-\rho) y]^2}
\ee
and use Proposition \ref{p:Fanal}(ii), to conclude from (\ref{fvarrred7}) that
$\frac{\partial}{\partial\rho}F(\rho)>0$. Hence (\ref{fvarred4}) reduces to
$F(\rho^*(p))$. 
\epr

The following proposition is the analogue of Proposition \ref{phtrcond}.

\bp{phtrcond*}
Fix $p<p_c$.\\
(i) If $\psi_{AB}(\bar x)=\psi_{AA}(\bar x)$ and $\psi_{BA}(\bar y)=\psi_{BB}(\bar y)$,
then $f=F(\rho^*(p))$.\\
(ii) If $\psi_{AB}(\bar x)>\psi_{AA}(\bar x)$ or $\psi_{BA}(\bar y)>\psi_{BB}(\bar y)$,
then $f>F(\rho^*(p))$.
\ep

\bpr
It follows from (\ref{fevar}) and (\ref{psiineq}) that $f$ is bounded below by the
right-hand side of (\ref{fvarred1}). The latter equals $F(\rho^*(p))$, as shown in 
the proof of Proposition \ref{p:fred}, and so
\be{flbgen}
f \geq F(\rho^*(p)).
\ee

\medskip\noindent
(i) Abbreviate $\theta_{kl}(a)=a\psi_{kl}(a)$. We know that $\theta_{AB}\geq
\theta_{AA}$ and $\theta_{BA}\geq\theta_{BB}$ (by (\ref{psiineq})), and 
that all four functions are concave (by Lemma \ref{psiprop}(i)). Since
$\theta_{AA}$ and $\theta_{BB}$ are both differentiable (by (\ref{ka}) and 
(\ref{ka1psi})), the assumption of equality at $\bar x$, respectively, $\bar y$ 
implies that $\theta_{AB}$ and $\theta_{BA}$ are differentiable at $\bar x$,
respectively, $\bar y$ and that the equality carries over to the derivatives. 
Thus, we have
\be{theeq1}
\theta_{AB}'(\bar x)=\theta_{AA}'(\bar x) \quad \mbox{ and } \quad
\theta_{BA}'(\bar y)=\theta_{BB}'(\bar y).
\ee
Fix $(\rho_{kl})\in\cA$. Abbreviate
\be{Nbardef}
\bar N = \rho_A \theta_{AA}(\bar x) + \rho_B \theta_{BB}(\bar y)
\quad \mbox{ and } \quad
\bar D = \rho_A \bar x + \rho_B \bar y.
\ee
The fact that $\bar x$ and $\bar y$ are the maximisers of
\be{varreds}
\sup_{a_{AA},a_{BB}\geq 2} 
\frac{\rho_A \theta_{AA}(a_{AA}) + \rho_B \theta_{BB}(a_{BB})}
{\rho_A a_{AA} + \rho_B a_{BB}}
\ee
implies that
\be{thetaeqs}
\theta_{AA}'(\bar x) = \theta_{BB}'(\bar y) = \frac{\bar N}{\bar D}.
\ee
Hence, all four derivatives in (\ref{theeq1}) are equal to $\bar N/\bar D$.
Next, abbreviate  
\be{Ndef}
N = \sum_{kl} \rho_{kl} \theta_{kl}(a_{kl}) \quad \mbox{ and } \quad
D = \sum_{kl} \rho_{kl} a_{kl}.
\ee
By the concavity of $a\mapsto\theta_{kl}(a)$ (recall Lemma \ref{psiprop}(i)),
we have that for all $(a_{kl})$,
\be{Nestin}
\begin{aligned}
N &\leq \bar N + \rho_{AA}(a_{AA}-\bar x)\theta_{AA}'(\bar x)
+ \rho_{AB}(a_{AB}-\bar x)\theta_{AB}'(\bar x)\\
&\qquad + \rho_{BA}(a_{BA}-\bar y)\theta_{BA}'(\bar y)
+ \rho_{BB}(a_{BB}-\bar y)\theta_{BB}'(\bar y)\\
&= \bar N + V\,\frac{\bar N}{\bar D},
\end{aligned}
\ee
where
\be{Vdef}
V = \rho_{AA}(a_{AA}-\bar x) + \rho_{AB}(a_{AB}-\bar x)
+ \rho_{BA}(a_{BA}-\bar y) + \rho_{BB}(a_{BB}-\bar y).
\ee
Moreover,
\be{Desteq}
D = \bar D + V.
\ee
Combining (\ref{Nestin}) and (\ref{Desteq}), we obtain
\be{NDestin}
\frac{N}{D} \leq \frac{\bar N + V\,\frac{\bar N}{\bar D}}{\bar D + V}
= \frac{\bar N}{\bar D}.
\ee
Thus we have proved that 
\be{finerhofix}
\sup_{(a_{kl})\in\cA} \frac{N}{D} \leq \frac{\bar N}{\bar D} = F(\rho_A).
\ee
Optimising over $(\rho_{kl})\in\cR$, we arrive at $f \leq F(\rho^*(p))$. 
Together with (\ref{flbgen}) this gives the claim.

\medskip\noindent
(ii) Suppose that $\psi_{AB}(\bar x)>\psi_{AA}(\bar x)$ or $\psi_{BA}(\bar y)
>\psi_{BB}(\bar y)$. Then, by picking $a_{AA}=a_{AB}=\bar x$, $a_{BA}=a_{BB}
=\bar y$ and $\rho_A=\rho^*(p)$ in (\ref{fevar}), we get
\be{festlb}
f > \frac{\rho^*(p)\bar x\psi_{AA}(\bar x) + (1-\rho^*(p))\bar y\psi_{BB}(\bar y)}
{\rho^*(p)\bar x + (1-\rho^*(p))\bar y} 
= F(\rho^*(p)),
\ee
where we use that $0<\rho^*(p)<1$ for all $0<p<p_c$ (recall Fig.\ 6).
\epr

\subsubsection{Proof of Theorem \ref{phtrcurve*}}
\label{S4.2.2}

The following proposition is the analogue of Proposition \ref{phtrest}.

\bp{p:phtrest*}
Fix $p<p_c$.\\
(i) $\partial\cD$ lies on or below the supercritical curve $\alpha\mapsto
\beta_c(\alpha)$.\\
(ii) $\{(\alpha,\alpha)\colon\,\alpha\in[0,\alpha^*]\}\subset\partial\cD$,
where $\alpha^*$ is the same constant as for the supercritical curve, but 
lies below the diagonal elsewhere.\\
(iii) For every $p<p_c$ there exists an $\alpha^*(p) \in (0,\infty)$ such that
the intersection of $\cD$ with the lower half of $\CONE$ is the line segment 
$\{(\beta+\alpha^*(p),\beta)\colon\,\beta\in [-\frac12 \alpha^*(p),0]\}$.\\
(iv) $\cD$ is continuous and each line $\{(\beta+C,\beta)\colon\,\beta\geq 0\}$, 
$C\in [-1,1)$, intersects $\partial\cD$ in at most one point.\\
(v) As $p \da 0$, the part of $\cD$ not containing the diagonal segment
of (ii) converges to the mirror image of the analytic continuation of
the supercritical curve outside $\CONE$. In particular, $\lim_{p\da 0}\alpha^*(p)
=\alpha_0$, with $\alpha_0$ the constant defined in the first line of 
{\rm (\ref{alpha01})}.\\
(vi) $\lim_{p\ua p_c}\alpha^*(p) = \alpha_1$, with $\alpha_1$ the constant
defined in the second line of {\rm (\ref{alpha01})}.
\ep

\bpr
(i) Abbreviate
\be{philb4}
G(\mu,a) = \frac12\left(\frac{\mu-1}{\mu}\right)\log\left(\frac{a}{a-2}\right)
+\frac{1}{\mu}\log[2(a-1)].
\ee
Then Proposition \ref{p:psiineqchar} can be rewritten as the statement that 
$\psi_{AB}(a)>\psi_{AA}(a)$ if and only if
\be{philbextext}
\phi^{\cI}(\mu) > \frac12\alpha + G(\mu,a) \quad \mbox{ for some } \mu\geq 1,
\ee
and similarly for $\psi_{BA}(a)>\psi_{BB}(a)$ with $\frac12\beta$ instead of
$\frac12\alpha$. To prove the claim, we must show that, for all $(\alpha,\beta)
\in\CONE$ and $p<p_c$, the following is true: For all $\mu\geq 1$,
\be{philb5*}
\phi^\cI(\mu) \leq \frac12\alpha+G(\mu,\bar x)
\quad \mbox{ and } \quad 
\phi^\cI(\mu) \leq \frac12\beta+G(\mu,\bar y)
\ee
imply
\be{philb5**}
\phi^\cI(\mu) \leq \frac12\alpha+G(\mu,a^*).
\ee
Indeed, by Theorem \ref{phtriden*}, Proposition \ref{p:fe} and Proposition 
\ref{p:phtrinfchar}, this yields that $(\alpha,\beta)\in\cD$ for $p<p_c$ 
implies $(\alpha,\beta)\in\cD$ for $p\geq p_c$.  
  
We first show that the first half of (\ref{philb5*}) is redundant. Indeed,
\be{philb7}
\begin{aligned}
&\left[\frac12\alpha+G(\mu,\bar x)\right]
-\left[\frac12\beta+G(\mu,\bar y)\right]\\
&\qquad = \frac{\alpha-\beta}{2} + [G(\mu,\bar x)-G(\mu,\bar y)]\\
&\qquad = -\frac12\log\left(\frac{\bar x(\bar y-2)}{\bar y(\bar x-2)}\right)
+\frac12\left(\frac{\mu-1}{\mu}\right)
\log\left(\frac{\bar x(\bar y-2)}{\bar y(\bar x-2)}\right) 
+\frac{1}{\mu}\log\left(\frac{\bar x-1}{\bar y-1}\right)\\
&\qquad = \frac{1}{2\mu}\left[\log\left(\frac{(\bar x-2)(\bar x-1)^2}{\bar x}\right)
-\log\left(\frac{(\bar y-2)(\bar y-1)^2}{\bar y}\right)\right]\\
&\qquad \geq 0,
\end{aligned}
\ee
where in the third line we use the second line of (\ref{xysol4}), and in the fifth
line we use that $\bar x \geq \bar y$ (recall Proposition \ref{p:Fanal}(i)). 

Thus, it remains to show that, for all $\mu\geq 1$,
\be{philb8}
\phi^\cI(\mu) \leq \frac12\beta+G(\mu,\bar y)
\quad \mbox{ implies } \quad \phi^\cI(\mu) \leq \frac12\alpha+G(\mu,a^*).
\ee
Indeed,
\be{philb9}
\begin{aligned}
&\left[\frac12\alpha+G(\mu,a^*)\right]
-\left[\frac12\beta+G(\mu,\bar y)\right]\\
&\qquad = \frac{\alpha-\beta}{2}+[G(\mu,a^*)-G(\mu,\bar y)]\\
&\qquad = -\frac12\log\left(\frac{\bar x(\bar y-2)}{\bar y(\bar x-2)}\right)
+\frac12\left(\frac{\mu-1}{\mu}\right)
\log\left(\frac{a^*(\bar y-2)}{\bar y(a^*-2)}\right)
+\frac{1}{\mu}\log\left(\frac{a^*-1}{\bar y-1}\right)\\
&\qquad = \frac12\left[\log\left(\frac{a^*}{a^*-2}\right)
-\log\left(\frac{\bar x}{\bar x-2}\right)\right]\\
&\qquad\qquad
+\frac{1}{2\mu}\left[\log\left(\frac{(a^*-2)(a^*-1)^2}{a^*}\right)
-\log\left(\frac{(\bar y-2)(\bar y-1)^2}{\bar y}\right)\right]\\
&\qquad \geq 0,
\end{aligned}
\ee
where in the third line we use the second line of (\ref{xysol4}), and in the 
fifth line we use that $\bar x \geq a^* \geq \bar y$ (recall Proposition 
\ref{p:Fanal}(i)) to get that both terms between square brackets are $\geq 0$. 

\medskip\noindent
{\bf Remark:} The redundancy of the first half of (\ref{philb5*}) shows that, 
as the critical curve is crossed from $\cD$ to $\cL$, localization 
occurs in the $BA$-blocks rather than in the $AB$-blocks. This is why the 
first criterion in Proposition \ref{phtrcond*}(i) is redundant (as was claimed 
in (\ref{DLdef*red})).

\medskip\noindent
(ii) If $\alpha=\beta$, then $\bar x=\bar y=a^*$ by Proposition \ref{p:Fanal}(ii).
Therefore the second criterion in Proposition \ref{phtrcond*}(i) reduces to 
$\psi_{BA}(a^*)=\psi_{BB}(a^*)$ (while the first criterion in Proposition 
\ref{phtrcond*}(i) is redundant). But if $\alpha=\beta$, then $\psi_{BA}=
\psi_{AB}$ and $\psi_{BB}=\psi_{AA}$. Hence, the criterion for delocalization
on the diagonal reads $\psi_{AB}(a^*)=\psi_{AA}(a^*)$, which is the same as 
the criterion for delocalization in the supercritical case (recall Proposition 
\ref{p:fe} and \ref{phtrcond}).

This also shows that $\partial\cD$ must leave the diagonal at the same point as 
the supercritical curve, i.e., at $(\alpha^*,\alpha^*)$. (Incidentally, note that 
$\alpha=\beta$ does not imply $\psi_{BA}=\psi_{BB}$  or $\psi_{AB}=\psi_{AA}$, 
because only matches of the polymer and the emulsion receive an energy.)

\medskip\noindent
(iii) Let $\alpha \geq 0$ and $\beta \leq 0$. By (\ref{phirelbds}), $\phi^{\cI}
(\mu;\alpha,\beta)=\frac12\alpha+\hat\kappa(\mu)$. It therefore follows from 
(\ref{DLdef*red}) and the line below (\ref{philbextext}) that $(\alpha,\beta)
\not\in\cD$ if and only if
\be{cutcond}
\exists\,\mu\geq 1\colon\quad \hat\kappa(\mu)
> G(\mu,\bar y) - \frac12 C,
\ee 
where we recall from Section \ref{S2.5} that $\bar y$ is a function of 
$C=\alpha-\beta$ and $\rho^*(p)$ only. Combining the second line of 
(\ref{xysol4}) with (\ref{philb4}), we have
\be{Grepalt}
\begin{aligned}
G(\mu,\bar y) - \frac12  C
&= \frac12\left(\frac{\mu-1}{\mu}\right)\log\left(\frac{\bar y}{\bar y-2}\right)
+ \frac{1}{\mu}\log[2(\bar y-1)] - \frac12 C\\
&= -\frac{1}{2\mu}\log\left(\frac{\bar y}{\bar y-2}\right)
+ \frac{1}{\mu}\log[2(\bar y-1)] +\frac12\log\left(\frac{\bar x}{\bar x-2}\right).
\end{aligned}
\ee
By Proposition \ref{p:Fanal}(iv), the right-hand side of (\ref{Grepalt}) is 
strictly decreasing in $C$ for fixed $\mu$ and $\rho^*(p)$. Hence there is 
a unique critical value $C^*$, which we call $\alpha^*(p)$, above which 
(\ref{cutcond}) holds. Since $\bar x \ua \infty$ and $\bar y \da 2$ as 
$C\to\infty$ by Proposition \ref{p:Fanal}(vii), the right-hand side tends to 
$-\infty$ as $-(1/2\mu)[\log(1/2(\bar y-2))+o(1)]+o(1)$. Since $\hat\kappa(\mu)
\geq 0$ for all $\mu\geq 1$, it follows that $C^*$ is finite.  

\medskip\noindent
(iv) The continuity of $\cD$ is immediate from Lemma \ref{psiprop}, Proposition
\ref{p:Fanal}(iv), Proposition \ref{Fanal} and Proposition \ref{phtrcond*}.  

To prove the remainder of the claim, we need the following.

\bl{monclaim}
$\beta \mapsto \phi^{\cI}(\mu;\beta+C,\beta)-\frac12\beta$ is non-decreasing on
$[0,\infty)$ for all $\mu\geq 1$ and $C\geq 0$.
\el

\bpr
The function $(\alpha,\beta)\mapsto\phi^{\cI}(\alpha,\beta;\mu)$ is convex on $\R^2$
for all $\mu\geq 1$, by an argument similar to that given in the proof of Theorem
\ref{feiden}(ii) in Section \ref{S3.1} (recall (\ref{feinf}--\ref{fesainf}) and 
(\ref{feinfscal})). Fix $\mu\geq 1$ and $C=\alpha-\beta\geq 0$. Abbreviate 
$\Delta(\beta)=\phi^{\cI}(\beta+C,\beta;\mu)-\frac12\beta$. Then $\beta\mapsto
\Delta(\beta)$ is convex. Moreover, by (\ref{phirelbds}), $\Delta(0)=\phi^{\cI}
(C,0;\mu)=\frac12 C+\hat\kappa(\mu)$ and $\Delta(\beta)\geq\frac12 (\beta+C)
+\hat\kappa(\mu)-\frac12\beta=\Delta(0)$ when $\beta\geq 0$. Hence 
$\beta\mapsto\Delta(\beta)$ is non-decreasing on $[0,\infty)$. 
\epr

\noindent
To prove the claim, use that $\bar x$ and $\bar y$ are functions of $C=\alpha-\beta$ 
only, repeat the same argument as in the proof of (iii), and use Lemma \ref{monclaim}. 

\medskip\noindent
(v) The limit $p \da 0$ corrsponds to $\rho \da 0$ (recall Fig.\ 6). By Proposition 
\ref{p:Fanal}(vi), $\bar y \ua a^*$ when $\alpha<\beta+\log 5$ and $\bar y \ua 
2/(1-e^{-(\alpha-\beta)})$ when $\alpha \geq \beta+\log 5$. In the first case, 
the criterion for delocalization reduces to
\be{critred}
\sup_{\mu\geq 1} \mu\left[\phi^{\cI}(\mu)-\frac12\beta-\frac12\log 5\right]
\leq \frac12\log\frac95.
\ee
Recalling Propositions \ref{p:fe} and \ref{p:phtrinfchar}, we see that this is 
precisely the criterion for delocalization in the superciritcal case but with 
$\alpha$ and $\beta$ interchanged. The second case is ruled out by the observation
that $\alpha_0 \leq \log 5$, as is immediate from the first line of (\ref{alpha01}).

\medskip\noindent
(vi) The limit $p \ua p_c$ corresponds to $\rho \ua 1$. Let $\beta=0$. Then, by 
Proposition \ref{p:Fanal}(v), $\bar y \da 10/(5-e^{-\alpha})$ as $\rho\ua 1$. 
Therefore, with the help of Proposition \ref{p:psiineqchar} and (\ref{phirelbds}), 
the criterion for delocalization when $\alpha \geq 0$ reduces to
\be{alphalimcrit}
\sup_{\mu\geq 1} \mu\left[\hat\kappa(\mu)-\frac12\log 5\right]
\leq \frac12\log\left[\frac{4e^{-\alpha}(5+e^{-\alpha})^2}
{5(5-e^{-\alpha})^2}\right].
\ee
This inequality holds if and only if $\alpha\leq\alpha_1$.
\epr

The redundancy of the first half of (\ref{philb5*}), which is strict when 
$\alpha>\beta$, shows that, as the critical curve is crossed from $\cD$ to 
$\cL$ off the diagonal, localization occurs in the $BA$-blocks rather than 
in the $AB$-blocks. 

\subsubsection{Critical lines on the diagonal}
\label{S4.2.3}

In this section we explain why in Figs.\ 9 and 10 the diagonal segment
$\{(\alpha,\alpha)\colon\,\alpha\in[-\alpha^*,\alpha^*]\}$ is a critical 
line. This is not obvious from our earlier considerations, because the
segment lies on the boundary of $\CONE$.

Take Fig.\ 9 ($p \geq p_c$). We have $f=\psi_{AA}(a^*)$ in the phase $\cD_A$ 
and $f=\psi_{AA}(\bar y)$ in the phase $\cD_{A+B}$, where $\bar y=\bar y(\beta,
\alpha;\rho^*(1-p))$ is the $y$-maximiser of (\ref{Fdef}) with $\alpha\leftrightarrow
\beta$ and $p\leftrightarrow 1-p$. Since $a^*$ is the unique maximiser of 
$\psi_{AA}$ (by Lemma \ref{l:ka}(iv) and (\ref{ka1psi})) and $\bar y \neq a^*$ 
when $\beta>\alpha$ (by Proposition \ref{p:Fanal}(i)), this shows that the 
free energy is non-analytic along the separation line between $\cD_A$ 
and $\cD_{A+B}$.

Take Fig.\ 10 ($1-p_c<p<p_c$). We have $f=F(\alpha,\beta;\rho^*(p))$ in the
lower half of the phase $\cD_{A+B}$ and $f=F(\beta,\alpha;\rho^*(1-p))$ in the 
upper half of the phase $\cD_{A+B}$. Similarly as in Theorem \ref{feiden}(iv),
there is the symmetry
\be{Fsymms}
F(\beta,\alpha;\rho^*(1-p)) = F(\alpha,\beta;1-\rho^*(1-p))
\ee 
(as is also evident from (\ref{Fdef}) and (\ref{xysol4})). However,
\be{rho*discr}
\rho^*(p) > 1-\rho^*(1-p) \qquad \forall\,p\in (0,1)
\ee
because the curve in Fig.\ 6 lies strictly above the diagonal. Since $\rho\mapsto
F(\alpha,\beta;\rho)$ is strictly increasing (as we saw below (\ref{fvarred8})), 
this shows that free energy is non-analytic along the separation line between 
the two halves of $\cD_{A+B}$.

It might be that in Fig.\ 10 also inside the two phases $\cL_{AB+BA}$ there is a 
critical line on the diagonal. We do not expect this, but we lack the tools to decide.

\subsection{Further observations about $\cL$ for $p<p_c$}
\label{S4.3}

We close by making some observations about the second subcritical curve,
lying inside $\cL$. By Proposition \ref{phtrcond*}(ii), the criterion for
$\cL$ is $\psi_{BA}(\bar y)>\psi_{BB}(\bar y)$, corresponding to
$BA$-localization. The phase $\cL$ splits further into two parts,
one where $AB$-localization does not occur and one where it does. The 
criterion for $AB$-localization reads 
\be{critL}
\psi_{AB}(\tilde x) > \psi_{AA}(\tilde x)
\ee
where we denote by $\tilde x$, $\tilde y$, $\tilde z$ the values of the 
minimisers $a_{AB}(=a_{AA})$, $a_{BA}$, $a_{BB}$, repectively, when in the 
variational expression for the free energy (\ref{fevar}) we replace 
$\psi_{AB}$ by $\psi_{AA}$. Indeed, this is in complete analogy with 
the argument in Section \ref{S4.2.1} identifying the first subcritical
curve as the one arising when in the variational expression for the 
free energy (\ref{fevar}) we replace $\psi_{AB}$ by $\psi_{AA}$ \emph{and}
$\psi_{BA}$ by $\psi_{BB}$. Unfortunately, whereas the latter reduction 
leads to a computable supremum (as shown in Section \ref{S4.2.1}), the
former reduction does not (because we have no closed form expression for
$\psi_{BA}$). Consequently, we have little information on $\tilde x$, 
$\tilde y$, $\tilde z$ (unlike for $\bar x$, $\bar y$), which is why 
(\ref{critL}) is hard to exploit. However, we can use (\ref{critL}) to 
obtain a lower bound on the second subcritical curve.

\bl{Lbd}
If $\alpha \geq 0$ and $\beta \leq \log (2-e^{-\alpha})$, then $\psi_{AB}
\equiv\psi_{AA}$ and hence {\rm (\ref{critL})} fails.
\el 

\bpr
From (\ref{Morita3}) we know that if $\beta\leq\log(2-e^{-\alpha})$, then
\be{phiub1}
\phi^\cI(\mu) \leq \frac12\alpha + \hat\kappa(\mu) 
\ee
and hence
\be{phiub2}
\sup_{\mu\geq 1} \mu\left[\phi^\cI(\mu)-\frac12\alpha
-\frac12\log\left(\frac{a}{a-2}\right)\right] \leq 
\sup_{\mu\geq 1} \mu\left[\hat\kappa(\mu)
-\frac12\log\left(\frac{a}{a-2}\right)\right].
\ee
But 
\be{phiub3}
\sup_{\mu\geq 1} \mu\left[\hat\kappa(\mu)
-\frac12\log\left(\frac{a}{a-2}\right)\right]
\leq \frac12\log\left[\frac{4(a-2)(a-1)^2}{a}\right],
\ee
as can be seen from Proposition \ref{p:psiineqchar}, because $\phi^\cI\equiv
\hat\kappa$ and $\psi_{AB}\equiv\psi_{AA}$ when $\alpha=\beta=0$. Combining 
(\ref{phiub2}) and (\ref{phiub3}) with Proposition \ref{p:psiineqchar}, we 
find that if $\beta\leq\log(2-e^{-\alpha})$, then $\psi_{AB}\equiv\psi_{AA}$. 
\epr

\noindent
Lemma \ref{Lbd} shows that the second critical curve is bounded below
by the curve $\alpha\mapsto\log(2-e^{-\alpha})$. 

The second subcritical curve splits off the diagonal at the same point 
$(\alpha^*,\alpha^*)$ as the first subcritical curve. Indeed, if $\alpha=\beta$, 
then $\psi_{AB}\equiv\psi_{BA}$, $\psi_{AA}\equiv\psi_{BB}$, and $\bar x=\bar y$. 
Therefore, on the diagonal the criteria for $AB$-localization ($\psi_{AB}(a^*)
=\psi_{AA}(a^*)$) and $BA$-localization ($\psi_{BA}(a^*)=\psi_{BB}(a^*)$) 
coincide.

We believe that the second critical curve has a finite horizontal asymptote,
as argued on physical grounds in Section \ref{S1.5extra}, and that it lies 
above the supercritical curve. We are unable to prove this.


\end{document}